\PassOptionsToPackage{stixpath=./latex/assets/fonts/}{walpuski-fonts}
\documentclass[
  annotations=show,
  font=stix2
]{walpuski-article}

\ifXeTeX
\fi

\NewDocumentCommand{\SpaceOfMetrics}{o}{\operatorname{\mathbf{Met}}\IfValueT{#1}{(#1)}}
\NewDocumentCommand{\SpaceOfClosedNonDegenerateMultiValuedOneForms}{O{} m}{\mu{}Z^1\IfBlankTF{#1}{}{_{#1}}(#2)}
\newcommand{\DeRhamCohomologyBundle}[1][\bullet]{\mathbf{H}_{\dR}^{#1}}
\newcommand{\GaussManinConnection}{\nabla} 
\newcommand{\PeriodSection}{\mathrm{Per}}
\NewDocumentCommand{\SpaceOfHarmonicNonDegenerateMultiValuedOneForms}{O{} m}{%
  \mu{}\IfBlankTF{#1}{}{_{#1}}\sH^1(#2)%
}
\newcommand{\RegularLocus}[1]{\operatorname{Reg}(#1)}
\NewDocumentCommand{\RamificationLocus}{O{} m}{\operatorname{Ram}\IfBlankTF{#1}{}{_{#1}}(#2)}
\newcommand{\BranchingLocus}[1]{\operatorname{Br}(#1)}
\NewDocumentCommand{\ModuliSpaceOfSpecialLagrangians}{o}{\sM_{\mathrm{SL}\IfValueT{#1}{,#1}}}
\newcommand{\DeformationNormalConeModuliSpaceOfSpecialLagrangians}{\sD_{\mathrm{SL}}}
\newcommand{\Current}[1]{\llbracket #1 \rrbracket}

\newcommand{\SpaceOfLagrangianImmersions}[2]{\operatorname{\mathbf{LagImm}}(#1,#2)}
\NewDocumentCommand{\SpaceOfLagrangianSubmanifolds}{O{} m}{\operatorname{\mathbf{LagSub}}\IfBlankTF{#1}{}{_{#1}}(#2)}
\newcommand{\SpaceOfParameterisedClosedNonDegenerateMultiValuedOneForms}[2]{\operatorname{\mathbf{LagImm}}^*(#1,T^*#2)}
\newcommand{\SetOfDiffeomorphismTypes}{\fL}
\newcommand{\ExactFormBundle}[1][\bullet]{\sB^{#1}}
\newcommand{\DefiningSection}{\eta}

\NewDocumentCommand{\PositiveLocus}{o}{\sU^+\IfValueT{#1}{_{#1}}}

\newcommand{\HodgeTypeOperator}[2]{H_{#1,#2}}

\newcommand{\iu}{\mathrm{i}}

\newcommand{\SubellipticLaplacian}{\square}
\NewDocumentCommand{\EnergyForm}{m o o}{\sE_{#1}\IfValueT{#2}{\paren{#2\IfValueT{#3}{,#3}}}}

\NewDocumentCommand{\dist}{o m m}{%
  \operatorname{dist}\IfValueT{#1}{_{#1}}\paren{#2,#3}%
}

\title[
  plaintext={Harmonic multi-valued 1–forms and deformations of branched covers of special Lagrangians}
]{
  Harmonic multi-valued $1$--forms and\\
  deformations of branched covers of special Lagrangians
}
\addauthor[
  inst={1},
  email={sqhe@amss.ac.cn},
  surname={He},
  onclick={https://siqihemath.github.io/Homepage/},
  orcid={0009-0003-4292-4762}
]{Siqi He}
\addauthor[
  inst={2},
  email={jacek.rzemieniecki@hu-berlin.de},
  surname={Rzemieniecki},
  onclick={https://jacek-rzemieniecki.net/},
  orcid={0009-0008-2195-8150}
]{Jacek Rzemieniecki}
\addauthor[
  inst={2},
  email={thomas.walpuski@hu-berlin.de},
  surname={Walpuski},
  onclick={https://walpu.ski/},
  orcid={0000-0002-2705-3423}
]{Thomas Walpuski}
\addaffiliation[
  department={Morningside Center of Mathematics},
  city={Beijing},
  country={China},
  countrycode={CN}
]{Academy of Mathematics and Systems Science, Chinese Academy of Sciences}
\addaffiliation[
  department={Institut für Mathematik},
  city={Berlin},
  country={Germany},
  countrycode={DE}
]{Humboldt-Universität zu Berlin}
\license{CC-BY-NC-SA-4.0}

\AtEveryBibitem{%
  \ifboolexpr{
    test {\iffieldequalstr{entrykey}{Mohsen2021:DeformationGroupoid}}
    or test {\iffieldequalstr{entrykey}{Taubes2013}}
  }
    {\emergencystretch=1em
     \renewbibmacro*{finentry}{\finentry\par}}
    {}%
}
\date{2026-11-24}

\begin{document}
\maketitle

\begin{abstract}
  This article introduces the notion of a non-degenerate multi-valued $1$--form on a manifold $X$ in terms of immersed submanifolds of $T^*X$ and
  develops a deformation theory extending Donaldson's work on harmonic non-degenerate $\ZModTwoZ$ $1$--forms.
  More precisely, the universal moduli space of harmonic non-degenerate multi-valued $1$--forms is shown to be a tame Fréchet orbifold whose projection onto the space of Riemannian metrics is a submersion and whose fibres are locally parametrised by their periods.
  The proof relies on Hörmander's subelliptic estimates and the Nash--Moser submersion theorem.
  As an application of this analytic framework,
  harmonic multi-valued $1$--forms on a special Lagrangian are shown to control deformations of its branched covers,
  extending work of He for multiplicity two.
\end{abstract}


\section{Introduction}
\label{Sec_Introduction}

His analysis of the failure of compactness for flat $\SL_2(\C)$--connections led \citet{Taubes2012} to introduce harmonic $\ZModTwoZ$ $1$--forms, that is,
harmonic two-valued $1$--forms with average zero.
Subsequently, they have appeared in related problems in gauge theory \cite{Taubes2013,Taubes2022:NahmPoleKapustinWitten} and in the context of branched double cover degenerations of special Lagrangians \cite{He2022:DeformBranchedSLags}.

Harmonic $1$--forms are governed by topology via Hodge theory.
It is natural to ask to what extent this remains true for harmonic $\ZModTwoZ$ $1$--forms, in particular, whether harmonic $\ZModTwoZ$ $1$--forms persist under small perturbations of the metric.
\citet[Theorem 1.1]{Donaldson2021} has developed such a deformation theory for harmonic $\ZModTwoZ$ $1$--forms which are \emph{non-degenerate}, that is,
their branching locus is smooth and they vanish precisely to order $\frac12$ along it.
His approach treats the branching locus as a free boundary of codimension two.

This article extends \cite[Theorem~1.1]{Donaldson2021} to harmonic multi-valued
$1$--forms of arbitrary multiplicity.
The notion of non-degeneracy is formulated in terms of the regularity of the graph of the multi-valued $1$--form,
without prescribing local models of the branching behaviour.
The linearisation of the resulting deformation theory is governed by subelliptic Laplacians.

\medskip

Throughout this article,
let $X$ be an oriented manifold of dimension $d \geq 2$.
Denote
by $\pi \co T^*X \to X$ the projection map,
by $\lambda \in \Omega^1(T^*X)$ the \defined{tautological $1$--form}, and
by $\omega \coloneq \rd\lambda$ the \defined{canonical symplectic form}.
Denote by $\SpaceOfMetrics[X]$ the space of Riemannian metrics on $X$;
if $X$ is closed,
then it is a tame Fréchet manifold.

\begin{definition}
  \label{Def_NonDegenerateMultiValuedOneForm}
  A \defined{non-degenerate multi-valued $1$--form} on $X$ is a non-empty $d$--dimensional immersed submanifold $[\iota \co L \imm T^*X]$ such that
  $\pi \circ \iota$ is proper and,
  for every $x \in L$,
  \begin{equation*}
    \dim \ker T_x \paren{\pi \circ \iota} \leq 2.
  \end{equation*}
\end{definition}

Every $1$--form $\alpha \in \Omega^1(X)$ corresponds to a non-degenerate multi-valued $1$--form $[\alpha \co X \imm T^*X]$.
The above permits the projection map $\pi \circ \iota \co L \to X$ to exhibit branching, but also further undesirable phenomena, e.g., folds.
The latter are ruled out a posteriori by imposing harmonicity;
see \autoref{Sec_StructureOfProjectionMap}.

\begin{definition}
  \label{Def_NonDegenerateMultiValuedOneForm_Properties}
  Let $[\iota \co L \imm T^*X]$ be a non-degenerate multi-valued $1$--form on $X$.
  \begin{enumerate}
  \item
    \label{Def_NonDegenerateMultiValuedOneForm_Properties_Closed}
    $[\iota]$ is \defined{closed} if it is Lagrangian, that is,
    \begin{equation*}
      \iota^*\omega = 0.
    \end{equation*}
  \item
    \label{Def_NonDegenerateMultiValuedOneForm_Properties_Coclosed}
    For every $g \in \SpaceOfMetrics[X]$,
    define $\Lambda_g \in \Omega^{d-1}(T^*X)$ by
    \begin{equation*}
      \Lambda_g\paren{\xi;v_1,\ldots,v_{d-1}} \coloneq (*_g\xi)\paren{T_\xi \pi(v_1),\ldots,T_\xi \pi(v_{d-1})}
    \end{equation*}
    and set
    \begin{equation*}
      \Omega_g \coloneq \rd\Lambda_g \in \Omega^d(T^*X).
    \end{equation*}
    $[\iota]$ is \defined{co-closed} with respect to $g$ if
    \begin{equation*}
      \iota^*\Omega_g = 0.
    \end{equation*}
  \item
    \label{Def_NonDegenerateMultiValuedOneForm_Properties_Harmonic}
    $[\iota]$ is \defined{harmonic} with respect to $g$ if it is closed and co-closed with respect to $g$.
  \end{enumerate}
\end{definition}

A $1$--form $\alpha \in \Omega^1(X)$ is harmonic in the usual sense if
and only if the corresponding non-degenerate multi-valued $1$--form
$[\alpha]$ is harmonic in the sense of
\autoref{Def_NonDegenerateMultiValuedOneForm_Properties},
because $\alpha^*\lambda = \alpha$ and $\alpha^*\Lambda_g = *_g\alpha$.
If $X$ is connected,
every harmonic non-degenerate multi-valued $1$--form
$[\iota\colon L \looparrowright T^*X]$ has a \defined{multiplicity}
$k = k[\iota] \coloneq \abs{\deg(\pi \circ \iota)} \in \N$ and
defines a harmonic $k$--valued $1$--form
$\alpha = \alpha([\iota]) \in W_\loc^{1,2}\Gamma\paren{X,\Sym^k T^*X}$ in the
sense of \citet{Almgren2000:BigRegularityPaper};
see \autoref{Sec_RelationHarmonicW12MultiValued1Forms}.
If $k = 2$ and $[\iota]$ is invariant under $\set{\pm 1} \circlearrowright T^*X$,
then it defines a harmonic non-degenerate $\ZModTwoZ$ $1$--form.
The converse also holds;
see \autoref{Sec_RelationHarmonicZModTwoZ1Forms}.

If $X$ is closed,
then the set of non-degenerate multi-valued closed $1$--forms on $X$,
denoted by $\SpaceOfClosedNonDegenerateMultiValuedOneForms{X}$,
is endowed with the structure of a tame Fréchet orbifold.%
\footnote{%
  Throughout this article,
  the prefix $\mu$ indicates the non-degenerate multi-valued analogue of the corresponding classical concept, e.g.,
  $Z^1(L) \coloneq \ker \paren{\rd \co \Omega^1(L) \to \Omega^2(L)}$.
}
There is a finite-rank vector bundle $\DeRhamCohomologyBundle[1]$ over $\SpaceOfClosedNonDegenerateMultiValuedOneForms{X}$ whose fibre over $[\iota \co L \imm T^*X]$ is isomorphic to $\rH_{\dR}^1(L)$.
Moreover, $\DeRhamCohomologyBundle[1]$ carries a \defined{Gauß--Manin connection} $\GaussManinConnection$ and admits a \defined{period section} $\PeriodSection \in \Gamma\paren{\SpaceOfClosedNonDegenerateMultiValuedOneForms{X},\DeRhamCohomologyBundle[1]}$ defined by $\PeriodSection([\iota]) \coloneq \sqparen[\big]{\iota,[\iota^*\lambda]}$;
see \autoref{Sec_SpaceOfNonDegenerateMultiValuedClosedOneForms}.

\begin{theorem}
  \label{Thm_UniversalSpaceOfNonDegenerateHarmonicOneForms}
  Assume that $X$ is closed.
  The \defined{universal moduli space of harmonic non-degenerate multi-valued $1$--forms}
  \begin{equation*}
    \SpaceOfHarmonicNonDegenerateMultiValuedOneForms{X}
    \coloneq
    \set[\big]{
      (g,[\iota]) \in \SpaceOfMetrics[X] \times \SpaceOfClosedNonDegenerateMultiValuedOneForms{X}
      :
      \iota^*\Omega_g = 0
    }
  \end{equation*}
  is a tame Fréchet suborbifold of $\SpaceOfMetrics[X] \times \SpaceOfClosedNonDegenerateMultiValuedOneForms{X}$ and
  the projection map
  \begin{equation*}
    \pr_1 \co \SpaceOfHarmonicNonDegenerateMultiValuedOneForms{X} \to \SpaceOfMetrics[X]
  \end{equation*}
  is a submersion\footnote{%
    Here and throughout this article,
    this means that $\pr_1$ is locally represented by a projection $F \times G \to F$ with $F$ and $G$ tame Fréchet spaces;
    see \cites[p.~2]{Gloeckner2016:Submersions}[Definition~1.49(c)]{Schmeding2023:InfiniteDimensionalDifferentialGeometry}.
  };
  moreover,
  the covariant derivative of the period section induces an isomorphism
  \begin{equation*}
    \GaussManinConnection\PeriodSection
    \co
    T\SpaceOfHarmonicNonDegenerateMultiValuedOneForms{X}/\SpaceOfMetrics[X]
    \iso
    \pr_2^*\DeRhamCohomologyBundle[1].
  \end{equation*}
\end{theorem}

\begin{remark}
  \label{Rmk_WhyTameFrechet}
  The linearisation of the defining equation $\iota^*\Omega_g = 0$ is governed by a second-order operator.
  This operator is elliptic on the \defined{regular locus}
  \begin{equation*}
    \RegularLocus{\iota}
    \coloneq
    \set[\big]{
      x \in L : T_x(\pi\circ\iota) ~\text{is invertible}
    },
  \end{equation*}
  but degenerates at the \defined{ramification locus} $\RamificationLocus{\iota} \coloneq L \setminus \RegularLocus{\iota}$---unless $d=2$.
  Nevertheless,
  it turns out to satisfy Hörmander's condition.
  Therefore,
  by Hörmander's subelliptic estimates \cite{Hormander1967:Hypoelliptic},
  the inverse does not gain two derivatives, but only $\epsilon > 0$.
  This loss of regularity prevents a (straightforward) application of the submersion theorem for Banach orbifolds, and
  justifies the above formulation in the tame Fréchet category.

  A loss of regularity also arises in \cite{Donaldson2021} and the related articles \cite{Parker2023:Deformation,HeParkerWalpuski2026} as a consequence of treating the branching locus as a free boundary of codimension two.
\end{remark}

\begin{remark}
  \label{Rmk_LocallyCyclic}
  For a harmonic non-degenerate $k$--valued $1$--form $[\iota]$ with $k \geq 3$,
  neither the ramification locus nor the \defined{branching locus} $\BranchingLocus{\iota} \coloneq \pi \circ \iota\paren{\RamificationLocus{\iota}}$ need be smooth.
  Indeed, consider $X \coloneq \C^2$ with its Euclidean metric.
  Identifying $T^*\C^2 \cong \C^2 \times \C^2$,
  for every $m \in \N$ and $\epsilon \in \R$,
  \begin{align*}
    \iota_{m, \epsilon}\co L \coloneq \C^2 &\to T^*\C^2\\
    (u, v) &\mapsto \paren[\big]{ u^3-3\epsilon uv^m, v; u, \tfrac32 \epsilon m u^2v^{m-1} }
  \end{align*}
  is a harmonic non-degenerate $3$--valued $1$--form with
  \begin{align*}
    \RamificationLocus{\iota_{m, \epsilon}}
    &=
    \set[\big]{ (u,v) \in L : u^2=\epsilon v^m }
    \qand \\
    \BranchingLocus{\iota_{m, \epsilon}}
    &=
    \set[\big]{ (z,w) \in X : z^2 = 4\epsilon^3 w^{3m} }.
  \end{align*}
  For $\epsilon = 0$,
  the projection is the standard cyclic covering $(u,v) \mapsto (u^3,v)$.
  Singular branching and ramification loci arise from arbitrarily small deformations of this model.
  Nevertheless,
  the \defined{simple ramification locus} is smooth by \autoref{Prop_StructureOfProjectionMap}~\autoref{Prop_StructureOfProjectionMap_SimpleRamification}.
\end{remark}

\medskip

The methods underpinning \autoref{Thm_UniversalSpaceOfNonDegenerateHarmonicOneForms} also enable the following application to special Lagrangians. Let $(Z,\omega,J,\theta)$ be a Calabi--Yau $d$--fold and
$[\check\iota \co \Lambda \imm Z]$ a closed immersed special Lagrangian submanifold,
that is, $\check\iota^*\omega = 0$ and $\check\iota^* \Im \theta = 0$.
\emph{Choose} an open neighbourhood $\im 0_\Lambda \subseteq U \subseteq T^*\Lambda$ satisfying $[0,1]\cdot U \subseteq U$ and a symplectic immersion
\begin{equation*}
  \jmath \co U \imm Z
  \inlinetext{such that}
  \jmath \circ 0_\Lambda = \check\iota.
\end{equation*}
For every $t \in \R^\times$,
denote by $t \co T^*\Lambda \to T^*\Lambda$ multiplication by $t$,
set $U_t \coloneq t^{-1}(U)$, and
define $\Omega_t \in \Omega^d\paren{U_t}$ by
\begin{equation}
  \label{Eq_OmegaT}
  \Omega_t \coloneq -t^{-1} t^*\jmath^*\Im\theta.
\end{equation}
If $[\iota_t \co L \imm T^*\Lambda] \in \SpaceOfClosedNonDegenerateMultiValuedOneForms{\Lambda}$ satisfies $\im \iota_t \subseteq U_t$,
then
$[\jmath \circ t \circ \iota_t \co L \imm Z]$ is special Lagrangian if and only if $\iota_t^*\Omega_t = 0$.
Therefore,
for every $t \in \R^\times$,
\begin{equation*}
  \ModuliSpaceOfSpecialLagrangians[t]
  \coloneq
  \set[\big]{
    [\iota] \in \SpaceOfClosedNonDegenerateMultiValuedOneForms{\Lambda}
    :
    \im \iota \subseteq U_t ~\text{and}~ \iota^*\Omega_t = 0
  }
\end{equation*}
maps to part of the \defined{moduli space of special Lagrangians} in $Z$,
which parametrises deformations of branched covers of $[\check\iota \co \Lambda \imm Z]$.
Evidently,
it parametrises \emph{exactly the same} deformations independent of $t \in \R^\times$.

A computation,
carried out in the proof of \autoref{Prop_OmegaT},
shows that
\begin{equation*}
  \lim_{t \to 0} \Omega_t = \Omega_g
\end{equation*}
with $g$ denoting the Riemannian metric on $\Lambda$ induced by $\omega$ and $J$.
With $U_0 \coloneq T^*\Lambda$ and $\Omega_0 \coloneq \Omega_g$,
the definition of $\ModuliSpaceOfSpecialLagrangians[t]$ extends to $t = 0$ and produces
\begin{equation*}
  \ModuliSpaceOfSpecialLagrangians[0] = \SpaceOfHarmonicNonDegenerateMultiValuedOneForms{\Lambda,g}.
\end{equation*}
Here $\SpaceOfHarmonicNonDegenerateMultiValuedOneForms{\Lambda,g}$ denotes the fibre of $\SpaceOfHarmonicNonDegenerateMultiValuedOneForms{\Lambda}$ over $g$.
\emph{Formally},
this can be regarded as the central fibre of a deformation to the normal cone;
see \cites[§2.2.1]{KashiwaraSchapira1985:MicrolocalStudyOfSheaves}[§1.1]{Mohsen2021:DeformationGroupoid}.

\begin{theorem}
  \label{Thm_BranchedCoversOfSpecialLagrangians}
  The total space
  \begin{equation*}
    \DeformationNormalConeModuliSpaceOfSpecialLagrangians
    \coloneq
    \coprod_{t \in \R} \ModuliSpaceOfSpecialLagrangians[t]
     =
    \set[\big]{
      (t,[\iota])
      \in \R \times \SpaceOfClosedNonDegenerateMultiValuedOneForms{\Lambda}
      :
      \im \iota \subseteq U_t
      ~\text{and}~
      \iota^*\Omega_t =0
    }
  \end{equation*}
  is a suborbifold of $\R \times \SpaceOfClosedNonDegenerateMultiValuedOneForms{\Lambda}$ and the projection map $\pr_1 \co \DeformationNormalConeModuliSpaceOfSpecialLagrangians \to \R$ is a submersion;
  moreover,
  the covariant derivative of the period section induces an isomorphism
  \begin{equation*}
    \GaussManinConnection\PeriodSection
    \co
    T \DeformationNormalConeModuliSpaceOfSpecialLagrangians/\R
    \iso
    \pr_2^*\DeRhamCohomologyBundle[1].
  \end{equation*}
\end{theorem}

As a consequence,
$\SpaceOfHarmonicNonDegenerateMultiValuedOneForms{\Lambda,g}$ controls nearby deformations of branched covers of $[\check\iota]$.
In particular, \autoref{Thm_BranchedCoversOfSpecialLagrangians} immediately implies the following extension of \cite[Theorem 1]{He2022:DeformBranchedSLags}.

\begin{cor}
  \label{Cor_BranchedCoversOfSpecialLagrangians}
  Suppose that $\Lambda$ is connected.
  For every $[\iota \co L \imm T^*\Lambda] \in \SpaceOfHarmonicNonDegenerateMultiValuedOneForms{\Lambda,g}$ with multiplicity $k \coloneq k[\iota] \in \N$,
  there is a $1$--parameter family $([\iota_t])_{0 < \abs{t} < T}$ of immersed special Lagrangians such that
  \begin{equation*}
    \lim_{t \to 0} \paren{\iota_t}_\sharp\Current{L} = k \cdot \check\iota_\sharp\Current{\Lambda}
  \end{equation*}
  as integral currents.
\end{cor}

\paragraph{Acknowledgements}
SH wishes to thank Shaoyun Bai, Rafe Mazzeo and Greg Parker for their helpful comments and discussions.
JR thanks Thibault Langlais for insightful conversations on related problems, which led him to expect a result such as \autoref{Lem_SmoothAfterPullback} to hold.
TW thanks Aleksander Doan for telling him the time in Okinawa.

This material is based in part upon work carried out while SH and TW were in residence at the Simons Laufer Mathematical Sciences Institute (previously known as MSRI) in Berkeley, California, during the Fall 2022 semester.
JR is supported by the Deutsche Forschungsgemeinschaft (DFG, German Research Foundation) under Germany's Excellence Strategy through the Berlin Mathematics Research Center MATH+ (EXC-2046/2, project ID: 390685689).

\paragraph{AI Declaration}
SH and TW have been investigating possible extensions of \cite[Theorem 1.1]{Donaldson2021} to harmonic multi-valued $1$--forms since their stay at SLMath in 2022.
These investigations attempted to work with definitions of harmonic non-degenerate multi-valued $1$--forms prescribing local models of the branching behaviour and extensions of the analytic theory developed in \cites{Donaldson2021,He2022:DeformBranchedSLags}.
This, however, turned out to be unsatisfactory in light of \autoref{Rmk_LocallyCyclic}.
The fundamental idea of this article,
to work on the graph and use Hörmander's subelliptic estimates,
was discovered independently by JR and AI (in the form of GPT-6 Astra).
The latter was fed with a manuscript detailing SH and TW's failed attempts.

AI has been used extensively in the preparation of this article
to review the literature,
develop definitions and proofs,
and proofread the manuscript.
The present article was written by the authors,
who take full responsibility for its content.


\section{The orbifold of closed non-degenerate multi-valued \texorpdfstring{$1$}{1}--forms}
\label{Sec_SpaceOfNonDegenerateMultiValuedClosedOneForms}

This section sets up the infinite-dimensional differential geometric framework within which the proofs of \autoref{Thm_UniversalSpaceOfNonDegenerateHarmonicOneForms} and \autoref{Thm_BranchedCoversOfSpecialLagrangians} take place.
This requires the theory of tame Fréchet manifolds,
to which \cite{Hamilton1982:NashMoser} remains an excellent introduction.

\begin{definition}
  Let $(Z,\omega)$ be a symplectic manifold of dimension $2d$ and
  $L$ a closed manifold of dimension $d$.
  \begin{enumerate}
  \item
    Denote by $\SpaceOfLagrangianImmersions{L}{Z} \subseteq C^\infty(L,Z)$ the subset of Lagrangian immersions $\iota \co L \imm Z$.
  \item
    $\Diff(L)$ acts on the right of $\SpaceOfLagrangianImmersions{L}{Z}$ by reparametrisation: $(\iota,\phi) \mapsto \iota \circ \phi$.
    Denote the quotient by
    \begin{equation*}
      \SpaceOfLagrangianSubmanifolds[L]{Z}
      \coloneq
      \SpaceOfLagrangianImmersions{L}{Z}/\Diff(L).
    \end{equation*}
  \item
    Denote by
    $\SpaceOfClosedNonDegenerateMultiValuedOneForms[L]{X}$ the subset of those
    $[\iota \co L \imm T^*X] \in \SpaceOfLagrangianSubmanifolds[L]{T^*X}$ that satisfy the condition in \autoref{Def_NonDegenerateMultiValuedOneForm},
    and by $\SpaceOfParameterisedClosedNonDegenerateMultiValuedOneForms{L}{X}$ its preimage in $\SpaceOfLagrangianImmersions{L}{T^*X}$.
  \end{enumerate}
  The \defined{set of closed non-degenerate multi-valued $1$--forms} on $X$ is the disjoint union
  \begin{equation*}
    \SpaceOfClosedNonDegenerateMultiValuedOneForms{X}
    \coloneq
    \coprod_{L \in \SetOfDiffeomorphismTypes} \SpaceOfClosedNonDegenerateMultiValuedOneForms[L]{X}.
  \end{equation*}
  Here $\SetOfDiffeomorphismTypes$ is a set of closed manifolds of dimension $d$ which contains precisely one representative of each diffeomorphism type.
  Of course, the precise choice of $\SetOfDiffeomorphismTypes$ is immaterial.
\end{definition}

For every closed manifold $L$,
$C^\infty(L,Z)$ is endowed with the structure of a tame Fréchet manifold \cite[Part II, Corollary 2.3.2]{Hamilton1982:NashMoser}, and $\Diff(L)$ is a tame Fréchet Lie group \cite[Part II, Theorem 2.3.5]{Hamilton1982:NashMoser}.
The following endows $\SpaceOfLagrangianSubmanifolds[L]{Z}$ with the structure of a tame Fréchet orbifold.

\begin{prop}[{cf.~\citet{CerveraMascaroMichor1991,MichorMumford2005,Lee2009:WeightedSubmanifolds}}]
  \label{Prop_SpaceOfLagrangianSubmanifolds}
  Let $(Z,\omega)$ be a symplectic manifold of dimension $2d$.
  For every closed manifold $L$ of dimension $d$,
  the following hold:
  \begin{enumerate}
  \item
    \label{Prop_SpaceOfLagrangianSubmanifolds_Parametrised}
    $\SpaceOfLagrangianImmersions{L}{Z}$ is a split tame Fréchet submanifold of $C^\infty(L,Z)$.
  \item
    \label{Prop_SpaceOfLagrangianSubmanifolds_Unparametrised}
    There is a unique tame Fréchet orbifold structure on $\SpaceOfLagrangianSubmanifolds[L]{Z}$ with respect to which the projection map $\SpaceOfLagrangianImmersions{L}{Z} \to \SpaceOfLagrangianSubmanifolds[L]{Z}$ is a $\Diff(L)$--principal orbibundle.
  \item
    \label{Prop_SpaceOfLagrangianSubmanifolds_TangentBundle}
    The map
    \begin{align*}
      T\SpaceOfLagrangianImmersions{L}{Z} &\to Z^1(L) \coloneq \ker \paren{\rd \co \Omega^1(L) \to \Omega^2(L)} \\
      (\iota,v) &\mapsto \iota^*\paren{i_v\omega}
    \end{align*}
    induces an isomorphism
    \begin{equation*}
      T\SpaceOfLagrangianSubmanifolds[L]{Z}
      \iso
      \SpaceOfLagrangianImmersions{L}{Z} \times_{\Diff(L)} Z^1(L).
    \end{equation*}
  \end{enumerate}
\end{prop}

Evidently,
$\SpaceOfParameterisedClosedNonDegenerateMultiValuedOneForms{L}{X} \subseteq \SpaceOfLagrangianImmersions{L}{T^*X}$ is open and $\Diff(L)$--invariant;
moreover,
$\SpaceOfClosedNonDegenerateMultiValuedOneForms[L]{X} \subseteq \SpaceOfLagrangianSubmanifolds[L]{T^*X}$ is open.
Therefore,
the above endows $\SpaceOfClosedNonDegenerateMultiValuedOneForms{X}$ with the structure of a tame Fréchet orbifold and,
for every $L \in \SetOfDiffeomorphismTypes$,
$\SpaceOfParameterisedClosedNonDegenerateMultiValuedOneForms{L}{X} \to \SpaceOfClosedNonDegenerateMultiValuedOneForms[L]{X}$ is a $\Diff(L)$--principal orbibundle.

\begin{proof}[Proof of \autoref{Prop_SpaceOfLagrangianSubmanifolds}]
  \begin{proofsteps}
  \item[\autoref{Prop_SpaceOfLagrangianSubmanifolds_Parametrised}]
    The following is a variation of the construction in \cite[§3, proof of Proposition 1]{Lee2009:WeightedSubmanifolds}.

    Let $\iota \in \SpaceOfLagrangianImmersions{L}{Z}$.
    By \cite[Lemma~3.1]{CerveraMascaroMichor1991}, its stabiliser $G_\iota \coloneq \Stab{\iota}{\Diff(L)}$ is finite.
    Using \cite[Theorem~4.19]{BatesWeinstein1997},
    construct a $G_\iota$--invariant open neighbourhood $\im 0_L \subseteq U \subseteq T^*L$ and a $G_\iota$--invariant symplectic immersion
    \begin{equation*}
      \jmath \co U \imm Z
      \inlinetext{
        such that
      }
      \jmath \circ 0_L = \iota.
    \end{equation*}

    Choose a $G_\iota$--invariant open neighbourhood $0_L \in \sU \subseteq \Gamma\paren{L,U} \subseteq \Omega^1(L)$.
    By \cite[Part II, Corollary 2.2.7 and Theorems 2.3.3 and 2.3.5]{Hamilton1982:NashMoser},
    \begin{align*}
      \sU \times \Diff(L) &\to C^\infty\paren{L,U} \\
      (\alpha,\phi) &\mapsto \alpha \circ \phi
    \end{align*}
    is a tame open embedding.
    For every $(\alpha,\phi) \in \sU \times \Diff(L)$,
    $\alpha \circ \phi$ is a Lagrangian immersion if and only if $\alpha \in Z^1(L)$.
    Since the splitting $\Omega^1(L) = Z^1(L) \oplus \rd^*\Omega^2(L)$ with respect to an auxiliary Riemannian metric on $L$ is tame,
    the subset of those $(\alpha,\phi) \in \sU \times \Diff(L)$ for which $\alpha \circ \phi$ is Lagrangian is a split tame Fréchet submanifold.

    By \cite[§2.4]{MichorMumford2005},
    there is a $G_\iota$--invariant open neighbourhood $0_L \in \sV \subseteq \sU$ such that
    \begin{align*}
      \paren{\sV \times \Diff(L)}/G_\iota &\to C^\infty(L,Z) \\
      [\alpha,\phi] &\mapsto \jmath \circ \alpha \circ \phi
    \end{align*}
    is a tame open embedding.
    Since $\jmath$ is symplectic,
    $\jmath \circ \alpha \circ \phi$ is Lagrangian if and only if $\alpha \circ \phi$ is.
    This proves the assertion.
  \item[\autoref{Prop_SpaceOfLagrangianSubmanifolds_Unparametrised}]
    The construction in \autoref{Prop_SpaceOfLagrangianSubmanifolds_Parametrised} proves that the projection map
    \begin{equation*}
      \SpaceOfLagrangianImmersions{L}{Z} \to \SpaceOfLagrangianSubmanifolds[L]{Z}
    \end{equation*}
    is locally modelled on
    \begin{equation*}
      \paren[\big]{\paren{\sV \cap Z^1(L)} \times \Diff(L)}/G_\iota \to \paren{\sV \cap Z^1(L)}/G_\iota.
    \end{equation*}
    Together with the details supplied in \cites[Proof of Theorem~1.5]{CerveraMascaroMichor1991}[§2.4]{MichorMumford2005},
    this proves the assertion.
  \item[\autoref{Prop_SpaceOfLagrangianSubmanifolds_TangentBundle}]
    This is immediate from the above presentation of the projection map.
    \qedhere
  \end{proofsteps}
\end{proof}

The space of immersed Lagrangian submanifolds
\begin{equation*}
  \SpaceOfLagrangianSubmanifolds{Z}
  \coloneq
  \coprod_{L \in \SetOfDiffeomorphismTypes}
  \SpaceOfLagrangianSubmanifolds[L]{Z}
\end{equation*}
carries the following graded local system.

\begin{definition}
  \label{Def_DeRhamCohomologyBundle}
  The \defined{de Rham cohomology bundle} is the graded vector orbibundle
  \begin{equation*}
    \DeRhamCohomologyBundle
    \coloneq
    \coprod_{L \in \SetOfDiffeomorphismTypes}
    \SpaceOfLagrangianImmersions{L}{Z} \times_{\Diff(L)} \rH_{\dR}^\bullet(L)
    \to
    \SpaceOfLagrangianSubmanifolds{Z}.
  \end{equation*}
  Since $\Diff_0(L)$ acts trivially on $\rH_{\dR}^\bullet(L)$,
  $\DeRhamCohomologyBundle$ has discrete structure group $\Diff(L)/\Diff_0(L)$.
  The \defined{Gauß--Manin connection} $\GaussManinConnection$ is the corresponding flat covariant derivative on $\DeRhamCohomologyBundle$.
\end{definition}

For $Z = T^*X$,
there is the following additional structure.

\begin{definition}
  \label{Def_PeridodSection+BundleOfExactForms+DefiningSection}
  \begin{enumerate}
  \item
    The \defined{period section} $\PeriodSection \in \Gamma\paren{\SpaceOfClosedNonDegenerateMultiValuedOneForms{X},\DeRhamCohomologyBundle[1]}$ is defined by
    \begin{equation*}
      \PeriodSection([\iota]) \coloneq \sqparen[\big]{\iota,\sqparen{\iota^*\lambda}}.
    \end{equation*}
  \item
    The \defined{bundle of exact forms} is the graded tame Fréchet space orbibundle
    \begin{equation*}
      \ExactFormBundle
      \coloneq
      \coprod_{L \in \SetOfDiffeomorphismTypes}
      \SpaceOfParameterisedClosedNonDegenerateMultiValuedOneForms{L}{X} \times_{\Diff(L)} \rd\Omega^{\bullet-1}(L)
      \to
      \SpaceOfClosedNonDegenerateMultiValuedOneForms{X}.
    \end{equation*}
    Define $\DefiningSection \in \Gamma\paren{\SpaceOfMetrics[X] \times \SpaceOfClosedNonDegenerateMultiValuedOneForms{X},\pr_2^*\ExactFormBundle[d]}$ by
    \begin{equation*}
      \DefiningSection(g,[\iota]) \coloneq \paren{g,[\iota;\iota^*\Omega_g]}.
    \end{equation*}
  \end{enumerate}
\end{definition}

Evidently,
\begin{equation*}
  \SpaceOfHarmonicNonDegenerateMultiValuedOneForms{X} = \DefiningSection^{-1}(0).
\end{equation*}
The next five sections culminate in the proof of the following result.

\begin{proposition}
  \label{Prop_ExtendedCovariantDerivativeIsIsomorphism}
  There are
  an open neighbourhood $\SpaceOfHarmonicNonDegenerateMultiValuedOneForms{X} \subseteq \sU \subseteq \SpaceOfMetrics[X] \times \SpaceOfClosedNonDegenerateMultiValuedOneForms{X}$ and
  a smooth tame covariant derivative $\nabla$ on $\pr_2^*\ExactFormBundle[d]|_{\sU}$
  such that
  \begin{equation*}
    \paren{\nabla\DefiningSection ~ \GaussManinConnection\PeriodSection}
    \co \pr_2^*T\SpaceOfClosedNonDegenerateMultiValuedOneForms{X}|_{\sU}
    \to
    \pr_2^*\paren{\ExactFormBundle[d] \oplus \DeRhamCohomologyBundle[1]}|_{\sU}
  \end{equation*}
  is an isomorphism of tame Fréchet space orbibundles.
\end{proposition}

\begin{proof}[Proof of \autoref{Thm_UniversalSpaceOfNonDegenerateHarmonicOneForms} assuming \autoref{Prop_ExtendedCovariantDerivativeIsIsomorphism}]
  This follows from the Nash--Moser submersion theorem in the form of \cite[Theorem 4.7]{HeParkerWalpuski2026}.
\end{proof}


\section{The positive locus and the positive connection}
\label{Sec_PositiveLocusAndPositiveConnection}

This section constructs the covariant derivative that appears in \autoref{Prop_ExtendedCovariantDerivativeIsIsomorphism}.
The construction is a variation of an idea due to \textcites[§5.3]{Solomon2013:CalabiHomomorphism}.
It is based on the following construction of a $d$--form on $T^*X$ that induces a volume form on the domain of every harmonic non-degenerate multi-valued $1$--form.

\begin{definition}
  \label{Def_Theta}
  Let $(V,g)$ be an oriented Euclidean vector space of dimension $d$.
  \begin{enumerate}
  \item
    Denote its volume form by $\vol_g \in \Wedge^d V^*$ and define the \defined{symplectic form} $\omega \in \Wedge^2 (V \oplus V^*)^*$ by
    \begin{equation*}
      \omega((v,\xi),(w,\eta)) \coloneq \xi(w) - \eta(v).
    \end{equation*}
  \item
    Consider the Hamiltonian flow $(\Phi_t)_{t \in \R}$ of the kinetic energy $(v,\xi) \mapsto \frac12\abs{\xi}^2$,
    that is,
    \begin{equation*}
      \Phi_t(v,\xi) \coloneq (v + t\xi^\sharp,\xi).
    \end{equation*}
    For every $i \in \N_0$,
    define $\nu_{g,i} \in \Wedge^d (V \oplus V^*)^*$ by
    \begin{equation*}
      \sum_{i=0}^\infty t^i \nu_{g,i}
      \coloneq
      \Phi_t^*\pr_V^*\vol_g.
    \end{equation*}
  \item
    \label{Def_Theta_Fibrewise}
    Set
    \begin{equation*}
      \Omega_g \coloneq \nu_{g,1} \qandq
      \Theta_g \coloneq \nu_{g,0} - \nu_{g,2}.
    \end{equation*}
  \end{enumerate}
  For every $g \in \SpaceOfMetrics[X]$,
  decompose $TT^*X \iso \pi^*\paren{TX \oplus T^*X}$ using the Levi-Civita connection.
  With respect to this identification, $\Omega_g \in \Omega^d(T^*X)$,
  introduced in \autoref{Def_NonDegenerateMultiValuedOneForm_Properties}~\autoref{Def_NonDegenerateMultiValuedOneForm_Properties_Coclosed},
  agrees fibrewise with $\Omega_g$ defined in \autoref{Def_Theta_Fibrewise}.
  Define
  \begin{equation*}
    \Theta_g \in \Omega^d(T^*X)
  \end{equation*}
  by the corresponding fibrewise construction.
\end{definition}

\begin{prop}
  \label{Prop_Harmonic=>AlmostHarmonic}
  For every $(g,[\iota \co L \imm T^*X]) \in \SpaceOfHarmonicNonDegenerateMultiValuedOneForms{X}$,
  $\iota^*\Theta_g \in \Omega^d(L)$ is a volume form.
  \qed
\end{prop}

This is an immediate consequence of the following.

\begin{prop}
  \label{Prop_Theta}
  Assume the situation of \autoref{Def_Theta}.
  If $L \subseteq V \oplus V^*$ is Lagrangian and $\Omega_g|_L = 0$,
  then:
  \begin{enumerate}
  \item
    \label{Prop_Theta_Harmonic=>CorankNeq1}
    $\dim\paren{L \cap V^*} \neq 1$, and
  \item
    \label{Prop_Theta_Positive=>CorankLeq2}
    $\Theta_g|_L \neq 0$ if and only if $\dim\paren{L \cap V^*} \in \set{0,2}$.
  \end{enumerate}
\end{prop}

\begin{proof}
  Set $P \coloneq \pr_V L$.
  Since $L$ is Lagrangian,
  $L \cap V^* = \Ann P$ and
  there is a self-adjoint ${S} \in \End(P)$ such that
  \begin{align*}
    {G} \co V = P \oplus P^\perp &\to L \\
    (v,w) &\mapsto \paren{v, \paren{{S}v + w}^\flat}
  \end{align*}
  is an isomorphism.
  Since
  \begin{equation*}
    \pr_V \circ \Phi_t \circ G
    =
    \paren{\one_{P} +t{S}} \oplus t\one_{P^\perp},
  \end{equation*}
  it follows that
  \begin{align*}
    {G}^*\Phi_t^*\pr_V^*\vol_g
    &=
      t^\delta \det(\one_P + t{S}) \vol_g
    \qwithq
    \delta \coloneq \dim\paren{L \cap V^*}.
  \end{align*}
  The factor expands as
  \begin{equation*}
    t^\delta\det(\one_P+tS) =
    t^\delta + \tr St^{\delta+1} + \tfrac12\paren[\big]{\paren{\tr S}^2 - \abs{S}^2} t^{\delta+2} + \cdots
  \end{equation*}
  with $\abs{S} \coloneq \sqrt{\tr(S^2)}$ denoting the Hilbert--Schmidt norm.
  Therefore,
  by comparing coefficients,
  \begin{equation}
    \label{Eq_OmegaThetaViaS}
    \frac{{G}^*\Omega_g}{\vol_g}
    =
    \begin{cases}
      \tr S & \textnormal{if}~\delta = 0 \\
      1 & \textnormal{if}~\delta = 1 \\
      0 & \textnormal{if}~\delta \geq 2
    \end{cases}
    \qandq
    \frac{{G}^*\Theta_g}{\vol_g}
    =
    \begin{cases}
      1 + \frac12\paren[\big]{\abs{{S}}^2 - \paren{\tr S}^2} & \textnormal{if}~\delta = 0 \\
      -\tr S & \textnormal{if}~\delta = 1 \\
      -1 & \textnormal{if}~\delta = 2 \\
      0 & \textnormal{if}~\delta \geq 3.
    \end{cases}
  \end{equation}
  In particular,
  if $\Omega_g|_L = 0$,
  then $\delta \neq 1$,
  that is,
  \autoref{Prop_Theta_Harmonic=>CorankNeq1} holds.
  Moreover, \autoref{Prop_Theta_Positive=>CorankLeq2} follows by direct inspection.
\end{proof}

\begin{definition}
  \label{Def_PositiveLocus}
  The forms $\Theta_g$ determine the following open subset and covariant derivative:
  \begin{enumerate}
  \item
    The \defined{positive locus} is the open subset
    \begin{equation*}
      \PositiveLocus
      \coloneq
      \set[\big]{
        (g,[\iota]) \in \SpaceOfMetrics[X] \times \SpaceOfClosedNonDegenerateMultiValuedOneForms{X}
        :
        \iota^*\Theta_g ~\text{is a volume form}
      }.
    \end{equation*}
    It decomposes as $\PositiveLocus = \coprod_{L \in \SetOfDiffeomorphismTypes} \PositiveLocus[L]$.
  \item
    For every $L \in \SetOfDiffeomorphismTypes$, the \defined{positive connection} is the $\Diff(L)$--principal connection $A$ on the bundle
    \begin{equation*}
      \pr_2^*\SpaceOfParameterisedClosedNonDegenerateMultiValuedOneForms{L}{X}|_{\PositiveLocus}
      \to \PositiveLocus[L]
    \end{equation*}
    defined by
    \begin{equation*}
      i_{A(\dot g,v)} \iota^*\Theta_g
      \coloneq
      \iota^*\paren{i_v\Theta_g}.
    \end{equation*}
    Equip $\pr_2^*\ExactFormBundle|_{\PositiveLocus}$ with the covariant derivative $\nabla$ induced by the positive connections.
  \end{enumerate}
\end{definition}


\section{The structure of the covariant derivative}
\label{Sec_StructureOfCovariantDerivative}

This section analyses the structure of the partial covariant derivative $\nabla \DefiningSection$.

\begin{prop}
  \label{Prop_CovariantDerivativeOfDefiningSection}
  For every $(g,[\iota \co L \imm T^*X]) \in \PositiveLocus$ and $\alpha \in Z^1(L)$,
  \begin{equation*}
    \nabla\DefiningSection_{(g,[\iota])}([\iota,\alpha])
    =
    \sqparen{\iota,\rd \HodgeTypeOperator{g}{\iota}\alpha}
  \end{equation*}
  with $\HodgeTypeOperator{g}{\iota} \in \Gamma\paren{L,\Hom\paren{T^*L,\Wedge^{d-1} T^*L}}$ defined by
  \begin{equation*}
    \HodgeTypeOperator{g}{\iota}\alpha
    \coloneq
    \iota^*(i_v\Omega_g)
    -
    \tfrac{\iota^*\Omega_g}{\iota^*\Theta_g} \iota^*(i_v\Theta_g)
    \qwithq
    \iota^*\paren{i_v\omega} = \alpha.
  \end{equation*}
\end{prop}

\begin{remark}
  \label{Rmk_CovariantDerivativeOfDefiningSection}
  The use of the positive connection essentially amounts to regarding the defining equation of $\SpaceOfHarmonicNonDegenerateMultiValuedOneForms{X}$ as $\frac{\iota^*\Omega_g}{\iota^*\Theta_g} = 0$, and
  the above is essentially the quotient rule.
\end{remark}

\begin{proof}[Proof of \autoref{Prop_CovariantDerivativeOfDefiningSection}]
  Evidently, $\HodgeTypeOperator{g}{\iota}\alpha$ is independent of the choice of $v$.
  It is convenient to carry out the computation on $\SpaceOfParameterisedClosedNonDegenerateMultiValuedOneForms{L}{X}$.
  By \autoref{Prop_SpaceOfLagrangianSubmanifolds}~\autoref{Prop_SpaceOfLagrangianSubmanifolds_TangentBundle},
  $v \in \Gamma\paren{L,\iota^*TT^*X}$ is a suitable lift of $[\iota,\alpha]$.
  Its $A$--horizontal component is
  \begin{equation*}
    v - \iota_*Av.
  \end{equation*}
  Since $\Omega_g$ is closed,
  by Cartan's formula,
  the derivative of the map $\iota \mapsto \iota^*\Omega_g$ in the direction of $v - \iota_*Av$ is precisely
  \begin{equation*}
    \rd\iota^*\paren{i_v\Omega_g} - \rd i_{Av}\iota^*\Omega_g.
  \end{equation*}
  By \autoref{Def_PositiveLocus},
  \begin{equation*}
    i_{Av}\iota^*\Theta_g = \iota^*\paren{i_v\Theta_g}.
  \end{equation*}
  Therefore,
  \begin{equation*}
    i_{Av}\iota^*\Omega_g
    =
    \frac{\iota^*\Omega_g}{\iota^*\Theta_g} i_{Av}\iota^*\Theta_g
    =
    \frac{\iota^*\Omega_g}{\iota^*\Theta_g} \iota^*\paren{i_v\Theta_g}.
    \qedhere
  \end{equation*}
\end{proof}

A moment's thought,
or a glance at \autoref{Prop_PointwiseHodgeTypeOperator},
reveals that if $[\iota]$ is harmonic,
then,
on $\RegularLocus{\iota}$,
$\HodgeTypeOperator{g}{\iota}$ is the Hodge--$\ast$--operator with respect to $\paren{\pi \circ \iota}^*g$ and, therefore, $f \mapsto \rd \HodgeTypeOperator{g}{\iota}\rd f$ is essentially the Laplacian;
however, on $\RamificationLocus{\iota}$ it degenerates---unless $d = 2$.
The remainder of this section establishes two of the three properties of $\HodgeTypeOperator{g}{\iota}$ required for the analysis in \autoref{Sec_FamiliesOfSubellipticLaplacians}:
symmetry and non-negativity.

\begin{prop}
  \label{Prop_HodgeTypeOperator}
  For every $(g,[\iota \co L \imm T^*X]) \in \PositiveLocus$,
  the following hold:
  \begin{enumerate}
  \item
    \label{Prop_HodgeTypeOperator_Symmetric}
    $\HodgeTypeOperator{g}{\iota}$ is symmetric,
    that is,
    for every $\alpha,\beta \in \Omega^1(L)$,
    \begin{equation*}
      \alpha \wedge \HodgeTypeOperator{g}{\iota} \beta
      =
      \beta \wedge \HodgeTypeOperator{g}{\iota} \alpha.
    \end{equation*}
  \item
    \label{Prop_HodgeTypeOperator_NonNegativeDefinite}
    $\HodgeTypeOperator{g}{\iota}$ is non-negative definite,
    that is,
    for every $\alpha \in \Omega^1(L)$,
    \begin{equation*}
      \frac{\alpha \wedge \HodgeTypeOperator{g}{\iota} \alpha}{\iota^*\Theta_g}
      \geq
      0.
    \end{equation*}
  \end{enumerate}
\end{prop}

This is an immediate consequence of the following.

\begin{prop}
  \label{Prop_PointwiseHodgeTypeOperator}
  Assume the situation of \autoref{Def_Theta}.
  If $L \subseteq V \oplus V^*$ is Lagrangian and $\Theta_g|_L \neq 0$,
  then the linear map $\HodgeTypeOperator{g}{L} \co L^* \to \Wedge^{d-1} L^*$ defined by
  \begin{equation*}
    \HodgeTypeOperator{g}{L}\alpha \coloneq \left.\sqparen*{i_{(v,\xi)} \paren*{\Omega_g - \tfrac{\Omega_g|_L}{\Theta_g|_L} \Theta_g}}\right|_L
    \qwithq
    i_{(v,\xi)}\omega|_L = \alpha
  \end{equation*}
  satisfies the following:
  \begin{enumerate}
  \item
    \label{Prop_PointwiseHodgeTypeOperator_Symmetric}
    $\HodgeTypeOperator{g}{L}$ is symmetric,
    that is,
    for every $\alpha,\beta \in L^*$,
    \begin{equation*}
      \alpha \wedge \HodgeTypeOperator{g}{L} \beta = \beta \wedge \HodgeTypeOperator{g}{L} \alpha.
    \end{equation*}
  \item
    \label{Prop_PointwiseHodgeTypeOperator_Domination}
    $\HodgeTypeOperator{g}{L}$ is non-negative definite,
    that is,
    for every $\alpha \in L^*$,
    \begin{equation*}
      \frac{\alpha \wedge \HodgeTypeOperator{g}{L}\alpha}{\Theta_g|_L} \geq 0.
    \end{equation*}
  \item
    \label{Prop_PointwiseHodgeTypeOperator_Harmonic}
    If $\Omega_g|_L = 0$,
    then the following hold:
    \begin{enumerate}
    \item
      \label{Prop_PointwiseHodgeTypeOperator_Harmonic_1}
      If $L \cap V^* = 0$,
      then $\pr_V|_L$ identifies $\HodgeTypeOperator{g}{L}$ with $* \co V^* \to \Wedge^{d-1} V^*$.
    \item
      \label{Prop_PointwiseHodgeTypeOperator_Harmonic_2}
      If $\dim (L \cap V^*) = 2$,
      then,
      for every $\alpha \in L^*$,
      \begin{equation*}
        \alpha \wedge \HodgeTypeOperator{g}{L} \alpha = \abs[\big]{\alpha|_{L\cap V^*}}^2 \Theta_g|_L.
      \end{equation*}
    \end{enumerate}
  \end{enumerate}
\end{prop}

\begin{proof}
  \begin{proofsteps}%
  \item[\autoref{Prop_PointwiseHodgeTypeOperator_Symmetric}]
    Since $\pr_V^*\vol_g \wedge \omega = 0$ and $(\Phi_t)_{t\in\R}$ is a Hamiltonian flow,
    \begin{equation*}
      \Omega_g \wedge \omega = 0 \qandq \Theta_g \wedge \omega = 0.
    \end{equation*}
    Therefore,
    if $i_{(v,\xi)}\omega|_L = \alpha$ and $i_{(w,\eta)}\omega|_L = \beta$,
    then
    \begin{equation*}
      0
      =
      \left.\paren*{i_{(v,\xi)}i_{(w,\eta)} \sqparen*{\paren*{\Omega_g - \tfrac{\Omega_g|_L}{\Theta_g|_L} \Theta_g}\wedge\omega}}\right|_L
      =
      \alpha \wedge \HodgeTypeOperator{g}{L}\beta - \beta \wedge \HodgeTypeOperator{g}{L}\alpha.
    \end{equation*}

  \item[\autoref{Prop_PointwiseHodgeTypeOperator_Domination}]
    For $0 < t \ll 1$,
    $\Phi_t(L) \cap V^* = 0$ and $\Theta_g|_{\Phi_t(L)} \neq 0$.
    Therefore, by continuity, it suffices to assume that $L \cap V^* = 0$.

    Let ${S}$ and ${G} \co V \iso L$ be as in the proof of \autoref{Prop_Theta}.
    For every $\xi \in V^*$ and $\abs{t} \ll 1$,
    \begin{align*}
      G^*\paren{i_{(0,\xi)}\Phi_t^*\pr_V^*\vol_g}
      &=
        G^*\Phi_t^*\pr_V^* i_{t\xi^\sharp}\vol_g \\
      &=
        (\one+t{S})^*i_{t\xi^\sharp}\vol_g
        =
        t\det\paren{\one+t{S}} i_{\paren{\one+t{S}}^{-1}\xi^\sharp}\vol_g.
    \end{align*}
    Since
    \begin{equation*}
      t\det\paren{\one + t{S}}\paren{\one + t{S}}^{-1}
      =
      t + \paren{(\tr S) \one - {S}} t^2 + \cdots,
    \end{equation*}
    by comparing coefficients,
    \begin{equation*}
      G^*\paren{i_{(0,\xi)}\Omega_g} = *_g\xi
      \qandq
      G^*\paren{i_{(0,\xi)}\Theta_g} = -*_g\paren*{(\tr S)\one - {S}}^*\xi.
    \end{equation*}
    Therefore,
    \begin{equation}
      \label{Eq_1015}
      G^*\HodgeTypeOperator{g}{L}G_* \xi
      =
      *_g \sqparen*{\one+ \tfrac{\tr S}{D}\paren{(\tr S) \one - S}}^*\xi
      \qwithq
      D \coloneq 1 + \tfrac12\abs{{S}}^2 - \tfrac12\paren{\tr S}^2.
    \end{equation}
    By direct computation and completing the square,
    \begin{align*}
      \frac{G_*\xi \wedge \HodgeTypeOperator{g}{L} G_*\xi}{\Theta_g|_L}
      &=
        \frac{D\abs{\xi}^2 + \tr S\inner{\xi}{(\tr S)\xi - {S}^*\xi}}{D^2} \\
      &=
        \frac{
        \abs{\xi}^2
        +
        \frac12\paren{\abs{S}^2\abs{\xi}^2 - \abs{S^*\xi}^2}
        +
        \frac12\abs{(\tr S)\xi - {S}^*\xi}^2
        }{D^2}.
    \end{align*}
  \item[\autoref{Prop_PointwiseHodgeTypeOperator_Harmonic_1}]
    If $\Omega_g|_L = 0$,
    then $\tr S = 0$ and the assertion follows from \autoref{Eq_1015}.
  \item[\autoref{Prop_PointwiseHodgeTypeOperator_Harmonic_2}]
    Let $P$, ${S}$ and ${G} \co V = P \oplus P^\perp \iso L = G(P) \oplus (L\cap V^*)$ be as in the proof of \autoref{Prop_Theta}.
    For every $(v,\xi) \in V \oplus V^*$ and $0 < \abs{t} \ll 1$,
    \begin{align*}
      G^*\paren{i_{(v,\xi)} \Phi_t^*\pr_V^*\vol_g}
      &=
        \sqparen{\paren{\one_P + t{S}} \oplus t\one_{P^\perp}}^* i_{v+t\xi^\sharp}\vol_g \\
      &=
        t^2 \det\paren{\one_P + t{S}} i_{\paren{\one + tS}^{-1}\pr_{P}\paren{v + t\xi^\sharp} + t^{-1}\pr_{P^\perp}\paren{v + t\xi^\sharp}} \vol_g.
    \end{align*}
    By comparing coefficients,
    \begin{equation*}
      G^*\paren{i_{(v,\xi)} \Omega_g} = i_{\pr_{P^\perp}v}\vol_g.
    \end{equation*}
    It follows from \autoref{Eq_OmegaThetaViaS}
    that $\alpha \coloneq \paren{i_{(v,\xi)}\omega}|_L \in L^*$ satisfies
    \begin{align*}
      G^*\paren{\alpha \wedge \HodgeTypeOperator{g}{L} \alpha}
      &=
        G^*\alpha \wedge i_{\pr_{P^\perp}v}\vol_g
        =
        (G^*\alpha)(\pr_{P^\perp} v) \vol_g \\
      &=
        - \abs{\pr_{P^\perp} v}^2 \vol_g
        =
        \abs[\big]{\alpha|_{L\cap V^*}}^2 G^*\Theta_g.
      \qedhere
    \end{align*}
  \end{proofsteps}
\end{proof}


\section{The structure of the projection map}
\label{Sec_StructureOfProjectionMap}

The analysis in \autoref{Sec_FamiliesOfSubellipticLaplacians} requires a further ingredient:
Hörmander's condition.
The latter is introduced in \autoref{Sec_HormandersCondition} and established using the following discussion of the structure of the projection map $\pi \circ \iota$.

\begin{lemma}
  \label{Lem_StructureOfProjectionMap}
  Let $[\iota \co L \imm T^*X]$ be a harmonic non-degenerate multi-valued $1$--form.
  For every $x \in \RamificationLocus{\iota}$,
  there are charts $\phi \co U \to \tilde U \subseteq \R^{d-2} \times \C$ and $\psi \co V \to \tilde V \subseteq \R^{d-2} \times \C$ of $L$ and $X$
  satisfying $x \in U$, $\pi \circ \iota(x) \in V$ and $\phi(x) = \psi\paren{\pi \circ \iota(x)} = 0$, and
  $e_x \in \set{2,3,\ldots}$,
  such that $\sigma \in C^\infty(\tilde U;\C)$ defined by
  \begin{equation*}
    (t,\sigma(t,z))
    \coloneq
    \psi \circ \pi \circ \iota \circ \phi^{-1}(t,z)
  \end{equation*}
  satisfies
  \begin{equation*}
    \sigma(0,z) = \bar z^{e_x} + O\paren{\abs{z}^{e_x+1}}.
  \end{equation*}
\end{lemma}

\begin{proof}[Proof of \autoref{Lem_StructureOfProjectionMap}]
  The proof has five steps.

  \begin{step}
    Construction of the charts.
  \end{step}

  Set $P \coloneq \im T_x (\pi \circ \iota)$.
  Choose a chart $\psi \co V \to \tilde V \subseteq P \oplus P^\perp = T_{\pi \circ \iota(x)}X$ of $X$
  with $\pi \circ \iota(x) \in V$ and
  \begin{equation*}
    \psi \circ \pi\circ \iota(x) = 0
    \qandq
    T_{\pi\circ\iota(x)} \psi = \one.
  \end{equation*}
  This induces the chart
  $\hat\psi \co T^*V \to T^*\tilde V \subseteq P \oplus P^\perp \oplus P^* \oplus \paren{P^\perp}^*$ of $T^*X$ defined by
  \begin{equation*}
    \hat\psi (y,\xi) \coloneq \paren{\psi(y), \xi \circ \paren{T_y\psi}^{-1} - \iota(x)}.
  \end{equation*}
  Since $T_x\iota\paren{\ker T_x(\pi\circ\iota)} = \Ann P$,
  \begin{equation*}
    T _x \sqparen{\paren{\pr_P \oplus \pr_{(P^\perp)^*}} \circ \hat \psi \circ \iota}
    \co
    T_xL \to P \oplus \paren{P^\perp}^*
  \end{equation*}
  is an isomorphism.
  Therefore,
  by the inverse function theorem,
  there is a chart $\phi \co U \to \tilde U \subseteq P \oplus \paren{P^\perp}^*$ with $x \in U$ and $\tilde U$ contractible
  such that
  \begin{equation*}
    \phi = \paren{\pr_P \oplus \pr_{(P^\perp)^*}} \circ \hat \psi \circ \iota.
  \end{equation*}

  \begin{step}
    The local generating function.
  \end{step}

  By construction,
  there are $\sigma \co \tilde U \to P^\perp$ and $\tau \co \tilde U \to P^*$ such that
  \begin{equation*}
    \hat\psi \circ \iota \circ \phi^{-1}(t,z)
    =
    \paren{ t, \sigma(t,z); \tau(t,z), z };
  \end{equation*}
  moreover,
  \begin{equation*}
    \sigma(0,0) = 0, \quad
    T_{(0,0)}\sigma = 0, \qandq
    \tau(0,0) = 0.
  \end{equation*}
  Since $\iota$ is Lagrangian,
  \begin{equation*}
    \rd\paren{\inner{\sigma}{\rd z} - \inner{\tau}{\rd t}} = 0.
  \end{equation*}
  Since $\tilde U$ is contractible,
  there is a unique $S \in C^\infty\paren{\tilde U}$ such that
  \begin{equation*}
    S(0,0) = 0, \quad
    \sigma = \del_zS, \qandq
    \tau = -\del_tS;
  \end{equation*}
  moreover, by construction,
  \begin{equation*}
    \del_z^2 S(0,0) = 0 \qandq \del_t\del_z S(0,0) = 0.
  \end{equation*}
  $S$ is a local generating function for $\iota$;
  cf.~\cite[Chapter 5, §1.4]{ArnoldGivental1990:SymplecticGeometry}.

  \begin{step}
    An elliptic PDE for $S(t,\cdot)$.
  \end{step}

  Given $(t,z) \in \tilde U$ and $Q \in \Hom\paren{S^2 \paren{P^\perp}^*,\R}$,
  define $\sF(t,z;Q) \in \R$ by
  \begin{equation*}
    \sF(t,z;Q) \, \paren{\phi_*\iota^*\Theta_g}_{(t,z)}
    \coloneq
    \paren{\hat\psi_*\Omega_g}_{\paren{t,\del_zS;-\del_tS,z}}
    \circ
    \Wedge^d
    \begin{pmatrix}
      \one & 0 \\
      \del_t\del_z S & Q \\
      -\del_t^2S & -\del_z\del_t S \\
      0 & \one
    \end{pmatrix}.
  \end{equation*}
  Since $\iota^*\Omega_g = 0$,
  for every $(t,z) \in \tilde U$,
  \begin{equation*}
    \sF(t,z;\del_z^2S(t,z)) = 0.
  \end{equation*}
  It remains to analyse the nature of this PDE for $S(t,\cdot)$.

  $\sF(t,z;0) = 0$,
  because $\hat\psi_*\Omega_g$ only has components of type $(d,0)$ and $(d-1,1)$ with respect to the decomposition $\sqparen{P \oplus P^\perp} \oplus \sqparen{P^* \oplus \paren{P^\perp}^*}$.
  The latter is a consequence of the fact that $\Omega_g$ has type $(d-1,1)$ with respect to the decomposition $TT^*X = \pi^*TX \oplus \pi^*T^*X$.
  Therefore,
  $\sF(t,z;Q)$ factorises as
  \begin{equation*}
    \sF\paren{t,z;Q}
    =
    a\paren{t,z;Q} Q.
  \end{equation*}

  By the computation in the proof of \autoref{Prop_PointwiseHodgeTypeOperator}~\autoref{Prop_PointwiseHodgeTypeOperator_Harmonic_2},
  \begin{equation*}
    a\paren{0,0;0} Q
    =
    \del_Q\sF(0,0;0) Q
    =
    - \tr_gQ.
  \end{equation*}
  Consequently,
  possibly after shrinking $\tilde U$ around $0$,
  the PDE
  \begin{equation}
    \label{Eq_GeneratingFunctionPDE}
    a\paren{t,z;\del_z^2S(t,z)}\del_z^2 S(t,z) = 0
  \end{equation}
  is uniformly elliptic.

  \begin{step}
    Either the assertion holds or $S(0,\cdot)$ vanishes.
  \end{step}

  By strong unique continuation \cite{Aronszajn1957},
  either $S(0,\cdot)$ vanishes or its Taylor expansion at $0$ has a non-zero term.
  In the latter case,
  the leading-order term is a homogeneous harmonic polynomial of degree at least three.
  Hence,
  choosing an identification $(P^\perp)^* = \C$,
  it is of the form $\Re(c z^{e_x+1})$ with $e_x \geq 2$ and $c \in \C^\times$.
  Therefore,
  with $z = z_1 + \iu z_2$,
  \begin{equation*}
    \sigma(0,z)
    = \del_{z_1}S(0,z) + \iu\del_{z_2}S(0,z)
    = (e_x+1) \bar c \bar z^{e_x} + O\paren{\abs{z}^{e_x+1}}.
  \end{equation*}
  The factor $(e_x+1)\bar c$ can be absorbed by a suitable modification of $\psi$ and $\phi$.

  \begin{step}
    $S(0,\cdot)$ does not vanish.
  \end{step}

  Suppose not.
  Denote by $F$ the connected component of $\paren{\pi \circ \iota}^{-1}\paren{\pi \circ \iota(x)} \subseteq L$ containing $x$.
  By the preceding step,
  for every $y \in F$,
  the analogue of $S(0,\cdot)$ vanishes and $F$ is locally defined by $t = 0$.
  Therefore, $F$ is a two-dimensional submanifold.
  Since $\pi \circ \iota$ is proper, $F$ is compact.

  By construction,
  $\iota \co F \to T_{\pi \circ \iota(x)}^*X$.
  After relabelling, assume that $\abs{\iota}^2$ attains its maximum on $F$ at $x$.
  By construction,
  \begin{equation*}
    \iota \circ \phi^{-1}(0,z)
    =
    \iota(x) + \paren{\tau(0,z),z}.
  \end{equation*}

  Since $S(0,\cdot) = 0$ and $\tau = - \del_tS$,
  by differentiating \autoref{Eq_GeneratingFunctionPDE} in $t$,
  \begin{equation*}
    a(0,0;0) \del_z^2 \tau(0,0) = 0.
  \end{equation*}
  Therefore,
  by the maximum principle,
  \begin{align*}
    0
    &\leq
      a(0,0;0) \left.\del_z^2 \abs{\iota \circ \phi^{-1}(0,z)}^2\right|_{z = 0} \\
    &=
      2 \sum_{i,j=1}^2 a_{ij}(0,0;0)
      \inner{\del_{z_i} \paren{\iota \circ \phi^{-1}}(0,0)}{\del_{z_j} \paren{\iota \circ \phi^{-1}}(0,0)}
      = -4,
  \end{align*}
  a contradiction.
\end{proof}

The following result is not needed to establish Hörmander's condition,
but is of independent interest and will be used in \autoref{Sec_RelationHarmonicW12MultiValued1Forms} and \autoref{Sec_RelationHarmonicZModTwoZ1Forms}.

\begin{prop}
  \label{Prop_StructureOfProjectionMap}
  For every $(g,[\iota \co L \imm T^*X]) \in \SpaceOfHarmonicNonDegenerateMultiValuedOneForms{X}$,
  the following hold:
  \begin{enumerate}
  \item
    \label{Prop_StructureOfProjectionMap_NoFolds}
    For every $x \in L$,
    \begin{equation*}
      \dim \ker T_x(\pi \circ \iota) \in \set{0,2}.
    \end{equation*}
  \item
    \label{Prop_StructureOfProjectionMap_RegularLocusIsOpenDense}
    The regular locus $\RegularLocus{\iota}$ is open and dense.
  \item
    \label{Prop_StructureOfProjectionMap_RamificationIndex}
    For every $x \in L$,
    \begin{equation*}
      e_x \coloneq \deg_x(\pi \circ \iota),
    \end{equation*}
    the \defined{ramification index} of $\pi \circ \iota$ at $x$,
    defined with respect to the orientation of $L$ induced by $\iota^*\Theta_g$,
    equals $1$ if $x \in \RegularLocus{\iota}$,
    and agrees with the exponent in \autoref{Lem_StructureOfProjectionMap} if $x \in \RamificationLocus{\iota}$.
    If $X$ is connected, then, for every $y \in X$,
    \begin{equation*}
      \sum_{x \in (\pi \circ \iota)^{-1}(y)} e_x
      =
      \deg(\pi \circ \iota).
    \end{equation*}
  \item
    \label{Prop_StructureOfProjectionMap_SimpleRamification}
    The \defined{simple ramification locus}
    \begin{equation*}
      \RamificationLocus[2]{\iota} \coloneq \set[\big]{ x \in \RamificationLocus{\iota}
        :
        e_x = 2}
    \end{equation*}
    is open in $\RamificationLocus{\iota}$ and a submanifold of codimension two in $L$;
    moreover,
    the restriction $\pi \circ \iota |_{\RamificationLocus[2]{\iota}}$ is an immersion.
  \end{enumerate}
\end{prop}

\begin{proof}
  \begin{proofsteps}
  \item[\autoref{Prop_StructureOfProjectionMap_NoFolds}]
    This is immediate from \autoref{Prop_Theta}~\autoref{Prop_Theta_Harmonic=>CorankNeq1}.
  \item[\autoref{Prop_StructureOfProjectionMap_RegularLocusIsOpenDense}]
    Evidently,
    $\RegularLocus{\iota}$ is open.
    It is dense because,
    for every $x \in \RamificationLocus{\iota}$,
    in the charts of \autoref{Lem_StructureOfProjectionMap},
    $\del_z\sigma(0,z)$ is invertible for every sufficiently small $z \neq 0$;
    therefore, $x \in \overline{\RegularLocus{\iota}}$.
  \item[\autoref{Prop_StructureOfProjectionMap_RamificationIndex}]
    This is immediate from \autoref{Lem_StructureOfProjectionMap}.
  \item[\autoref{Prop_StructureOfProjectionMap_SimpleRamification}]
    Given $x \in \RamificationLocus[2]{\iota}$,
    choose charts as in \autoref{Lem_StructureOfProjectionMap},
    in which
    \begin{equation*}
      (t,\sigma(t,z))
      = \psi \circ \pi \circ \iota \circ \phi^{-1}(t,z)
      \qandq
      \sigma(0,z) = \bar z^2 + O\paren{\abs{z}^3}.
    \end{equation*}
    By \autoref{Prop_StructureOfProjectionMap_NoFolds},
    with $z = z_1 + \iu z_2$,
    $\phi^{-1}(t,z) \in \RamificationLocus{\iota}$ if and only if $\del_{z_1}\sigma(t,z) = 0$ if and only if $\del_{z_2}\sigma(t,z) = 0$.
    By inspection,
    \begin{equation*}
      \del_z\paren{\del_{z_1}\sigma}(0,0)\zeta = 2\bar\zeta.
    \end{equation*}
    Therefore,
    \begin{equation*}
      (0,0)
      \in
      \tilde W
      \coloneq
      \set[\big]{
        (t,z) \in \tilde U : \del_z\paren{\del_{z_1}\sigma}(t,z) ~\text{is invertible}
      }
      \qandq
      x \in W \coloneq \phi^{-1}(\tilde W).
    \end{equation*}
    By the regular value theorem,
    $\RamificationLocus{\iota} \cap W$ is a submanifold of codimension two.
    By \autoref{Lem_StructureOfProjectionMap},
    $\RamificationLocus{\iota} \cap W \subseteq \RamificationLocus[2]{\iota}$.
    Direct inspection shows that $\pi \circ \iota|_{\RamificationLocus{\iota} \cap W}$ is an immersion.
    \qedhere
  \end{proofsteps}
\end{proof}


\section{Hörmander's condition}
\label{Sec_HormandersCondition}

The following,
together with the two properties derived in \autoref{Prop_HodgeTypeOperator},
enables the analysis in \autoref{Sec_FamiliesOfSubellipticLaplacians}.

\begin{definition}
  \label{Def_HormandersCondition}
  Let $L$ be a manifold of dimension $d$ and $H \in \Gamma\paren{L,\Hom\paren{T^*L, \Wedge^{d-1}T^*L}}$.
  Define
  \begin{equation*}
    V_H \in \Gamma\paren{L,\Hom\paren{T^*L \otimes \Wedge^d TL,TL}}
  \end{equation*}
  by
  \begin{equation*}
    V_H(\alpha \otimes \tau) \coloneq i_{H\alpha} \tau,
    \inlinetext{and set}
    \sV_H \coloneq \im V_H \subseteq \Vect(L).
  \end{equation*}
  $H$ satisfies \defined{Hörmander's condition} if
  \begin{equation*}
    \overline{\sV}_H^{\Lie} = \Vect(L).
  \end{equation*}
  Here $\overline{\sV}_H^{\Lie}$ denotes the Lie-algebraic closure,
  that is, the smallest $C^\infty$--submodule which contains $\sV_H$ and is a Lie subalgebra of $\Vect(L)$.
\end{definition}

\begin{prop}
  \label{Prop_CoefficientModuleSatisfiesHormandersCondition}
  For every $(g,[\iota \co L \imm T^*X]) \in \SpaceOfHarmonicNonDegenerateMultiValuedOneForms{X}$,
  $\HodgeTypeOperator{g}{\iota}$ satisfies Hörmander's condition.
\end{prop}

\begin{proof}
  The proof has three steps.
  Choose a nowhere-vanishing $\tau \in \Gamma\paren{L,\Wedge^d TL}$ and set
  \begin{equation*}
    \rho \coloneq (\pi\circ \iota)^*\vol_g(\tau).
  \end{equation*}

  \begin{step}
    $\rho \, \Vect(L) \subseteq \sV_{\HodgeTypeOperator{g}{\iota}}$.
  \end{step}

  By \autoref{Prop_PointwiseHodgeTypeOperator}~\autoref{Prop_PointwiseHodgeTypeOperator_Harmonic_1},
  on the regular locus, for every $\xi \in \Omega^1(X)$,
  \begin{equation*}
    \HodgeTypeOperator{g}{\iota}(\pi\circ\iota)^*\xi
    =
    (\pi \circ \iota)^* *_g\xi.
  \end{equation*}
  Therefore, on the regular locus, for every $v \in \Vect(L)$,
  \begin{equation}
    \label{Eq_VAndVolume}
    V_{\HodgeTypeOperator{g}{\iota}}\paren{\sqparen{(\pi \circ \iota)^*g}(v,\cdot) \otimes \tau}
    =
    (\pi\circ \iota)^*\vol_g(\tau) v
    =
    \rho v.
  \end{equation}
  At critical points, $\rho = 0$ and
  $\sqparen{(\pi \circ \iota)^*g}(v,\cdot)$ annihilates $\ker T(\pi \circ \iota)$;
  therefore, by \autoref{Prop_PointwiseHodgeTypeOperator}~\autoref{Prop_PointwiseHodgeTypeOperator_Harmonic_2},
  both sides of \autoref{Eq_VAndVolume} vanish.

  \begin{step}
    For every $x \in L$, there are
    $v \in \sV_{\HodgeTypeOperator{g}{\iota}}$ and $k \in \N_0$
    such that $v^k\rho(x) \neq 0$.
  \end{step}

  At regular points, $\rho(x) \neq 0$;
  therefore, the assertion holds with $k = 0$.
  Suppose that $x \in \RamificationLocus{\iota}$.
  By \autoref{Prop_PointwiseHodgeTypeOperator}~\autoref{Prop_PointwiseHodgeTypeOperator_Harmonic_2},
  choose $\alpha \in \Omega^1(L)$ such that
  \begin{equation*}
    v \coloneq V_{\HodgeTypeOperator{g}{\iota}}(\alpha \otimes \tau)
    \qquad\text{satisfies}\qquad
    0 \neq v(x) \in \ker T_x(\pi \circ \iota).
  \end{equation*}
  By \autoref{Prop_HodgeTypeOperator}~\autoref{Prop_HodgeTypeOperator_Symmetric}
  and \autoref{Eq_VAndVolume},
  \begin{align*}
    \abs{(\pi \circ \iota)_*v}_g^2
    &=
      \sqparen{(\pi \circ \iota)^*g}(v,v) \\
    &=
      \alpha\paren{
      V_{\HodgeTypeOperator{g}{\iota}}\paren{
      \sqparen{(\pi \circ \iota)^*g}(v,\cdot) \otimes \tau
      }
      } \\
    &=
      \rho \alpha(v).
  \end{align*}

  Suppose that $v^k\rho(x) = 0$ for every $k \in \N_0$.
  Denote by $\gamma \co (-\epsilon,\epsilon) \to L$ an integral curve of $v$ with $\gamma(0) = x$.
  By construction,
  $\rho \circ \gamma$ vanishes to infinite order at $0$.
  Therefore,
  by the above,
  $\pi \circ \iota \circ \gamma$ is constant to infinite order at $0$.

  In the charts from \autoref{Lem_StructureOfProjectionMap},
  set
  \begin{equation*}
    (t(s),z(s)) \coloneq \phi \circ \gamma(s).
  \end{equation*}
  Both $t(s)$ and $\sigma(t(s),z(s))$ vanish to infinite order at $s = 0$,
  but
  $z'(0) \neq 0$,
  because $\gamma'(0) = v(x) \neq 0$ and $t'(0) = 0$.
  This, however, contradicts
  \begin{align*}
    \sigma(t(s),z(s))
    &=
      \sigma(0,z(s)) + O\paren{\abs{t(s)}} \\
    &=
      \overline{z'(0)}^{e_x} s^{e_x} + O\paren{\abs{s}^{e_x+1}}.
  \end{align*}

  \begin{step}
    $\overline{\sV}_{\HodgeTypeOperator{g}{\iota}}^{\Lie} = \Vect(L)$.
  \end{step}

  Let $x \in L$.
  Using the preceding step,
  choose $v \in \overline{\sV}_{\HodgeTypeOperator{g}{\iota}}^{\Lie}$ and
  let $k \in \N_0$ be the least integer
  such that
  $v^k\rho(x) \neq 0$.
  For every $w \in \Vect(L)$,
  \begin{equation*}
    \sL_v^k(\rho w)(x)
    =
    v^k\rho(x) \, w(x).
  \end{equation*}
  Since $\sL_v^k(\rho w) \in \overline{\sV}_{\HodgeTypeOperator{g}{\iota}}^{\Lie}$,
  $w(x)$ is arbitrary, and $v^k\rho(x) \neq 0$,
  it follows that $\overline{\sV}_{\HodgeTypeOperator{g}{\iota}}^{\Lie} = \Vect(L)$.
\end{proof}


\section{Families of subelliptic Laplacians}
\label{Sec_FamiliesOfSubellipticLaplacians}

This section proves the following result and, with its help, completes the proofs of \autoref{Prop_ExtendedCovariantDerivativeIsIsomorphism} and, therefore, \autoref{Thm_UniversalSpaceOfNonDegenerateHarmonicOneForms}.

\begin{theorem}
  \label{Thm_TameFamilyOfInversesForSubellipticLaplacians}
  Let $L$ be an oriented closed manifold of dimension $d$,
  ${\sP}$ be a tame Fréchet manifold,
  $p_0 \in {\sP}$, and
  $H \co \sP \to \Gamma\paren{L,\Hom\paren{T^*L,\Wedge^{d-1} T^*L}}$ a smooth tame map.
  If
  \begin{enumerate}
  \item
    for every $p \in {\sP}$,
    $H_p \coloneq H(p) \in \Gamma\paren{L,\Hom\paren{T^*L,\Wedge^{d-1} T^*L}}$ is \defined{symmetric} and \defined{non-negative definite},
    that is:
    \begin{enumerate}
    \item
      \label{Thm_TameFamilyOfInversesForSubellipticLaplacians_Symmetric}
      for every $\alpha, \beta \in \Omega^1(L)$,
      \begin{equation*}
        \alpha \wedge H_p \beta
        =
        \beta \wedge H_p \alpha,
      \end{equation*}
      and
    \item
      \label{Thm_TameFamilyOfInversesForSubellipticLaplacians_NonNegativeDefinite}
      for every $\alpha \in \Omega^1(L)$,
      \begin{equation*}
        \alpha \wedge H_p \alpha \geq 0,
      \end{equation*}
    \end{enumerate}
    and
  \item
    \label{Thm_TameFamilyOfInversesForSubellipticLaplacians_HormandersCondition}
    $H_{p_0}$ satisfies Hörmander's condition in the sense of \autoref{Def_HormandersCondition},
  \end{enumerate}
  then there is an open neighbourhood $p_0 \in \sU \subseteq \sP$
  such that,
  for every $p \in {\sU}$,
  \begin{equation*}
    \SubellipticLaplacian_p \coloneq \rd \circ H_p \co \rd\Omega^0(L) \to \rd\Omega^{d-1}(L)
  \end{equation*}
  is invertible and the inverses form a smooth tame map $\sU \times \rd\Omega^{d-1}(L) \to \rd\Omega^0(L)$.
\end{theorem}

\begin{proof}[Proof of \autoref{Prop_ExtendedCovariantDerivativeIsIsomorphism} assuming \autoref{Thm_TameFamilyOfInversesForSubellipticLaplacians}]
  By
  \autoref{Prop_CovariantDerivativeOfDefiningSection},
  \autoref{Prop_HodgeTypeOperator},
  \autoref{Prop_CoefficientModuleSatisfiesHormandersCondition}, and
  \autoref{Thm_TameFamilyOfInversesForSubellipticLaplacians},
  there is an open neighbourhood $\SpaceOfHarmonicNonDegenerateMultiValuedOneForms{X} \subseteq {\sU} \subseteq \PositiveLocus$ such that
  \begin{equation*}
    \nabla\DefiningSection
    \co
    \pr_2^*\ExactFormBundle[1]|_{\sU}
    \to
    \pr_2^*\ExactFormBundle[d]|_{\sU}
  \end{equation*}
  is an isomorphism of tame Fréchet space orbibundles.

  With respect to the isomorphism from \autoref{Prop_SpaceOfLagrangianSubmanifolds}~\autoref{Prop_SpaceOfLagrangianSubmanifolds_TangentBundle},
  $\GaussManinConnection\PeriodSection \co T\SpaceOfClosedNonDegenerateMultiValuedOneForms{X} \to \DeRhamCohomologyBundle[1]$
  is given by $\sqparen{\iota,\alpha} \mapsto \sqparen[\big]{\iota,[\alpha]}$.
  Therefore, the following commutative diagram
  \begin{equation*}
    \begin{tikzcd}
      \pr_2^*\ExactFormBundle[1]|_{\sU} \ar{r}{\nabla\DefiningSection} \ar[hook]{d}& \pr_2^*\ExactFormBundle[d]|_{\sU} \ar[hook]{d} \\
      \pr_2^*T\SpaceOfClosedNonDegenerateMultiValuedOneForms{X}|_{\sU} \ar{r}{\paren{\nabla\DefiningSection ~ \GaussManinConnection\PeriodSection}} \ar[two heads]{d}{\GaussManinConnection\PeriodSection} & \pr_2^*\paren{\ExactFormBundle[d] \oplus \DeRhamCohomologyBundle[1]} |_{\sU} \ar[two heads]{d} \\
      \pr_2^*\DeRhamCohomologyBundle[1]|_{\sU} \ar[equals]{r} & \pr_2^*\DeRhamCohomologyBundle[1]|_{\sU}
    \end{tikzcd}
  \end{equation*}
  has split exact columns.
  Consequently,
  $\paren{\nabla\DefiningSection ~ \GaussManinConnection\PeriodSection}$ is an isomorphism of tame Fréchet space orbibundles.
\end{proof}

The proof of \autoref{Thm_TameFamilyOfInversesForSubellipticLaplacians} and its preparation occupy the remainder of this section.
It is convenient to choose a Riemannian metric $g$ on $L$.
This determines the subspace
\begin{equation*}
  C_\perp^\infty(L)
  \coloneq
  \set[\big]{
    f \in C^\infty(L)
    :
    \inner{f}{c}_{L^2} = 0
    ~\text{for every locally constant}~c \in C^\infty(L)
  }
  \subseteq
  C^\infty(L).
\end{equation*}
The maps
\begin{equation*}
  \begin{split}
    C_\perp^\infty(L) &\to \rd\Omega^0(L) \\
    f &\mapsto \rd f
  \end{split}
  \qandq
  \begin{split}
    C_\perp^\infty(L) &\to \rd\Omega^{d-1}(L) \\
    f &\mapsto -f \, \vol_g
  \end{split}
\end{equation*}
are tame isomorphisms.
Through these isomorphisms $\SubellipticLaplacian_p \co \rd\Omega^0(L) \to \rd\Omega^{d-1}(L)$ corresponds to the \defined{subelliptic Laplacian} $\SubellipticLaplacian_p \co C_\perp^\infty(L) \to C_\perp^\infty(L)$ defined by
\begin{equation*}
  \SubellipticLaplacian_p f \coloneq - \frac{\rd H_p \rd f}{\vol_g}.
\end{equation*}
It is convenient to extend the latter to $\SubellipticLaplacian_p \co C^\infty(L) \to C_\perp^\infty(L)$ and consider the symmetric bilinear form $\EnergyForm{p} \co C^\infty(L) \times C^\infty(L) \to \R$ defined by
\begin{equation*}
  \EnergyForm{p}[f_1][f_2]
  \coloneq
  \int_L \rd f_1 \wedge H_p \rd f_2
  =
  \inner{f_1}{\SubellipticLaplacian_p f_2}_{L^2}.
\end{equation*}

The following subelliptic estimate and corresponding Neumann--Poincaré inequality lie at the heart of the invertibility assertion in \autoref{Thm_TameFamilyOfInversesForSubellipticLaplacians}.

\begin{proposition}
  \label{Prop_SubellipticEstimate+NeumannPoincareInequality}
  There are
  an open neighbourhood $p_0 \in \sV \subseteq \sP$ and $\epsilon > 0$
  such that,
  for every $p \in \sV$,
  the following hold:
  \begin{enumerate}
  \item
    \label{Prop_SubellipticEstimate+NeumannPoincareInequality_S}
    For every $f \in C^\infty(L)$,
    \begin{equation*}
      \norm{f}_{H^\epsilon}^2
      \lesssim
      \EnergyForm{p}[f][f] + \norm{f}_{L^2}^2.
    \end{equation*}
  \item
    \label{Prop_SubellipticEstimate+NeumannPoincareInequality_N}
    For every $f \in C_\perp^\infty(L)$,
    \begin{equation*}
      \norm{f}_{H^\epsilon}^2
      \lesssim
      \EnergyForm{p}[f][f].
    \end{equation*}
  \end{enumerate}
\end{proposition}

\begin{proof}
  \begin{proofsteps}
  \item[\autoref{Prop_SubellipticEstimate+NeumannPoincareInequality_S}]
    Choose $\alpha_1, \dots, \alpha_N \in \Omega^1(L)$ spanning $T^*L$.
    For every $p \in {\sP}$,
    the vector fields
    \begin{equation*}
      v_{p,i} \coloneq V_{H_p} (\alpha_i \otimes \vol_g^{-1}) \in \Vect(L)
      \quad
      (i \in \set{1,\ldots,N})
    \end{equation*}
    generate $\sV_{H_p}$,
    defined in \autoref{Def_HormandersCondition}.

    For every $p \in \sP$, $k \in \N$, and $I \in \set{1,\ldots,N}^k$,
    define $v_{p,I} \in \Vect(L)$ recursively by
    \begin{equation*}
      v_{p,(i)} \coloneq v_{p,i}
      \qandq
      v_{p,(i_1,\ldots,i_k)} \coloneq \sqparen{v_{p,i_1},v_{p,(i_2,\ldots,i_k)}}.
    \end{equation*}
    Since $H_{p_0}$ satisfies Hörmander's condition,
    there is an $m \in \N$ such that,
    for $p = p_0$,
    \begin{equation*}
      v_{p,I} \quad
      (k \in \set{1,\ldots,m} \qandq I \in \set{1,\ldots,N}^k)
    \end{equation*}
    span $TL$.
    Evidently, there is an open neighbourhood $p_0 \in \sV' \subseteq \sP$ such that the same holds for every $p \in \sV'$.
    In fact,
    given $r \in \N$,
    $\sV'$ can be chosen such that
    \begin{equation}
      \label{Eq_UniformityInput}
      \begin{split}
        \inf_{p \in \sV'} \min_{\substack{\xi \in T^*L \\ \abs{\xi} = 1}}
        \sum_{k=1}^m \sum_{I \in \set{1,\ldots,N}^k} \abs{\xi\paren{v_{p,I}}}^2
        &> 0
          \qand \\
        \sup_{p \in \sV'} \max_{i \in \set{1,\ldots,N}} \norm{v_{p,i}}_{C^r}
        &<
          \infty.
      \end{split}
    \end{equation}

    Consequently,
    by Hörmander's subelliptic estimates \cite[Theorem 4.3]{Hormander1967:Hypoelliptic},
    there is an $\epsilon > 0$
    such that,
    for every $p \in \sV'$ and $f \in C^\infty(L)$,
    \begin{equation*}
      \norm{f}^2_{H^\epsilon}
      \lesssim
      \sum_{i=1}^{N} \norm{v_{p,i}f}_{L^2}^2 + \norm{f}^2_{L^2}.
    \end{equation*}
    \citeauthor{Hormander1967:Hypoelliptic} does not state a uniform bound,
    but,
    for sufficiently large $r \in \N$,
    uniformity in $p \in \sV'$ follows from \autoref{Eq_UniformityInput} by inspection of the proof of \cite[Theorem 4.3]{Hormander1967:Hypoelliptic}.
    (A modern exposition of both Hörmander's proof and Kohn's proof \cite{Kohn1973:Hypoellipticity} can be found in \cite[§2 and §3]{Bramanti2020:ProofOfHormandersTheorem}.)

    By the Cauchy--Schwarz inequality,
    applied to the non-negative definite symmetric bilinear form $(\alpha,\beta) \mapsto \frac{\alpha \wedge H_p \beta}{\vol_g}$,
    for every $p \in {\sP}$, $i \in \set{1,\ldots,N}$, and $f \in C^\infty(L)$,
    \begin{equation*}
      \abs{v_{p,i} f}^2
      =
      \paren*{\frac{\rd f \wedge H_p \alpha_i}{\vol_g}}^2
      \leq
      \frac{\alpha_i \wedge H_p \alpha_i}{\vol_g}
      \frac{\rd f \wedge H_p \rd f}{\vol_g}.
    \end{equation*}
    Therefore,
    there is an open neighbourhood $p_0 \in \sV'' \subseteq \sP$
    such that,
    for every $p \in \sV''$ and $f \in C^\infty(L)$,
    \begin{equation}
      \label{Eq_EnergyEstimate}
      \sum_{i=1}^N \norm{v_{p, i}f}^2_{L^2}
      \lesssim
      \EnergyForm{p}[f][f].
    \end{equation}
    This proves \autoref{Prop_SubellipticEstimate+NeumannPoincareInequality_S} for every open neighbourhood $p_0 \in \sV \subseteq \sV' \cap \sV''$.
  \item[\autoref{Prop_SubellipticEstimate+NeumannPoincareInequality_N}]
    If it is not possible to choose $\sV$ such that \autoref{Prop_SubellipticEstimate+NeumannPoincareInequality_N} holds,
    then there exist
    a sequence $(p_n) \in \sP^\N$ converging to $p_0 \in \sP$ and
    a sequence $(f_n) \in C_\perp^\infty(L)^\N$ such that
    \begin{equation*}
      \norm{f_n}_{L^2} = 1
      \qandq
      \lim_{n \to \infty} \EnergyForm{p_n}[f_n][f_n] = 0.
    \end{equation*}
    By \autoref{Prop_SubellipticEstimate+NeumannPoincareInequality_S} and since $H^\epsilon(L) \incl L^2(L)$ is compact,
    after passing to a subsequence,
    $(f_n)$ converges to a non-zero limit $f \in L^2(L)$ which is $L^2$--orthogonal to every locally constant $c \in C^\infty(L)$.
    Moreover, by \autoref{Eq_EnergyEstimate},
    for every $i \in \set{1,\ldots,N}$,
    \begin{equation*}
      v_{p_0,i} f = 0 \quad\text{in the sense of distributions.}
    \end{equation*}
    Since $H_{p_0}$ satisfies Hörmander's condition,
    for every $v \in \Vect(L)$,
    $vf = 0$ in the sense of distributions.
    It follows that $f$ is itself locally constant and, therefore, $f = 0$:
    a contradiction.
    \qedhere
  \end{proofsteps}
\end{proof}

\autoref{Prop_SubellipticEstimate+NeumannPoincareInequality} immediately implies that $\SubellipticLaplacian_p \co C_\perp^\infty(L) \to C_\perp^\infty(L)$ is injective.
Together with the Riesz representation theorem and subelliptic regularity, it also implies surjectivity; see, e.g., \cite[Proposition 1.9]{Baudoin2022:GeometricInequalities}.

The alternative argument presented below avoids appealing to subelliptic regularity and proceeds by elliptic regularisation \cite[§4.1 and §4.2]{KohnNirenberg1965}.
For this purpose,
for every $t \in [0,1]$,
set
\begin{equation*}
  \SubellipticLaplacian_{p,t} \coloneq \SubellipticLaplacian_p + t\Delta_g \co C^\infty(L) \to C_\perp^\infty(L)
\end{equation*}
and consider the symmetric bilinear form $\EnergyForm{p,t} \co C^\infty(L) \times C^\infty(L) \to \R$ defined by
\begin{equation*}
  \EnergyForm{p,t}[f_1][f_2]
  \coloneq
  \EnergyForm{p}[f_1][f_2] + t \inner{\nabla f_1}{\nabla f_2}_{L^2}
  =
  \inner{f_1}{\SubellipticLaplacian_{p,t} f_2}_{L^2}.
\end{equation*}
The following tame estimate is the remaining ingredient of the proof of \autoref{Thm_TameFamilyOfInversesForSubellipticLaplacians}.

\begin{proposition}
  \label{Prop_TameEstimate}
  There is an open neighbourhood $p_0 \in \sU \subseteq \sV$ such that,
  for every $p \in \sU$, $t \in [0,1]$, $k \in \N_0$, and $f \in C_\perp^\infty(L)$,
  \begin{equation*}
    \norm{f}_{H^{k+\epsilon}}
    \lesssim_k
    \norm{\SubellipticLaplacian_{p,t} f}_{H^k} + \norm{H_p}_{C^{k+2}}\norm{\SubellipticLaplacian_{p,t} f}_{H^1}.
  \end{equation*}
\end{proposition}

\begin{proof}
  The following is a variation of \cite[Proof of Proposition 2.27]{DonaldsonLehmann2024:Closed3FormInFiveDimensions}.

  By \autoref{Prop_SubellipticEstimate+NeumannPoincareInequality}~\autoref{Prop_SubellipticEstimate+NeumannPoincareInequality_N}
  and the Cauchy--Schwarz inequality,
  \begin{equation*}
    \norm{f}_{H^\epsilon}^2
    \lesssim \EnergyForm{p}[f][f]
    \leq \EnergyForm{p,t}[f][f]
    = \inner{f}{\SubellipticLaplacian_{p,t} f}_{L^2}
    \leq \norm{f}_{L^2}\norm{\SubellipticLaplacian_{p,t} f}_{L^2}.
  \end{equation*}
  Therefore,
  \begin{equation}
    \label{Eq_TameEstimate_K=0}
    \norm{f}_{H^\epsilon} \lesssim \norm{\SubellipticLaplacian_{p,t} f}_{L^2}.
  \end{equation}
  This proves the assertion for $k = 0$.

  The proof for $k \in \N$ combines this argument with what is essentially an elaborate commutator estimate.
  As in \cite[p.~420]{Hamilton1979:DeformationOfComplexStructuresII},
  for every $k \in \N$,
  choose $N = N(k) \in \N$ and $v_1, \ldots, v_N \in \Vect(L)$ such that,
  for every $f \in C^\infty(L)$,
  \begin{equation*}
    \norm{f}_{H^{k+\epsilon}}^2
    \leq
    \sum_{i=1}^{N}
    \sum_{j=0}^k
    \norm{v_i^j f}_{H^\epsilon}^2.
  \end{equation*}

  By \autoref{Prop_SubellipticEstimate+NeumannPoincareInequality}~\autoref{Prop_SubellipticEstimate+NeumannPoincareInequality_S},
  for every $p \in \sV$,
  \begin{equation*}
    \norm{f}_{H^{k+\epsilon}}^2
    \lesssim_k
    \sum_{i=1}^N
    \sum_{j=1}^k
    \EnergyForm{p,t}[v_i^jf][v_i^jf] + \norm{f}_{H^k}^2.
  \end{equation*}
  By \citeauthor{Hamilton1979:DeformationOfComplexStructuresII}'s uniform Kohn--Nirenberg lemma \cite[§4.5]{Hamilton1979:DeformationOfComplexStructuresII},
  there is an open neighbourhood $p_0 \in \sU \subseteq \sV$,
  independent of $k \in \N$,
  such that,
  for every $p \in \sU$, $t \in [0,1]$, $i \in \set{1,\ldots,N}$, and $j \in \set{1,\ldots,k}$,
  \begin{align*}
    \abs[\big]{\EnergyForm{p,t}[v_i^jf][v_i^jf] - (-1)^j\EnergyForm{p,t}[f][v_i^{2j}f]}
    &\lesssim_k
      \norm{f}_{H^j}^2 + \norm{H_p + t *_g}_{C^{j+2}}^2\norm{f}_{H^1}^2 \\
    &\lesssim_k
      \norm{f}_{H^k}^2 + \norm{H_p}_{C^{k+2}}^2\norm{f}_{H^1}^2.
  \end{align*}
  By integration by parts,
  \begin{align*}
    \abs[\big]{\EnergyForm{p,t}[f][v_i^{2j} f]}
    =
    \abs[\big]{\inner{\SubellipticLaplacian_{p,t} f}{v_i^{2j}f}_{L^2}}
    &=
      \abs[\big]{\inner{({v}_i^*)^j\SubellipticLaplacian_{p,t} f}{{v}_i^j f}_{L^2}} \\
    &\lesssim_k
      \norm{\SubellipticLaplacian_{p,t} f}_{H^j} \norm{f}_{H^j}
      \lesssim
      \norm{\SubellipticLaplacian_{p,t} f}_{H^k}^2 + \norm{f}_{H^k}^2.
  \end{align*}
  Therefore,
  \begin{equation*}
    \norm{f}_{H^{k+\epsilon}}^2
    \lesssim_k
    \norm{\SubellipticLaplacian_{p,t} f}_{H^k}^2
    + \norm{f}_{H^k}^2 + \norm{H_p}_{C^{k+2}}^2\norm{f}_{H^1}^2.
  \end{equation*}
  It remains to absorb the second term and control the third term.

  By Sobolev interpolation and \autoref{Eq_TameEstimate_K=0},
  for every $\delta > 0$,
  there is a $c(k,\delta) > 0$ such that
  \begin{equation*}
    \norm{f}_{H^k}
    \leq
    \delta \norm{f}_{H^{k+\epsilon}}
    + c(k,\delta) \norm{\SubellipticLaplacian_{p,t} f}_{L^2}.
  \end{equation*}
  For $\delta \ll 1$,
  rearranging absorbs the second term.

  The choice of $\sU$ in \cite[§4.5]{Hamilton1979:DeformationOfComplexStructuresII} entails that $\sup_{p \in \sU} \norm{H_p}_{C^3} < \infty$.
  Therefore,
  \begin{equation*}
    \norm{f}_{H^{1+\epsilon}}
    \lesssim
    \norm{\SubellipticLaplacian_{p,t} f}_{H^1}
    +
    \norm{f}_{H^1}.
  \end{equation*}
  By Sobolev interpolation, \autoref{Eq_TameEstimate_K=0}, and rearranging,
  \begin{equation*}
    \norm{f}_{H^1}
    \lesssim
    \norm{\SubellipticLaplacian_{p,t} f}_{H^1}.
  \end{equation*}
  Therefore, the third term is controlled.
  \qedhere
\end{proof}

\begin{proof}[Proof of \autoref{Thm_TameFamilyOfInversesForSubellipticLaplacians}]
  The problem has already been reduced to one about the family of subelliptic Laplacians $\paren[\big]{\SubellipticLaplacian_p \co C_\perp^\infty(L) \to C_\perp^\infty(L)}_{p \in \sU}$.
  The remainder of the proof consists of two steps.

  \begin{step}
    For every $p \in \sU$,
    $\SubellipticLaplacian_p \co C_\perp^\infty(L) \to C_\perp^\infty(L)$ is a tame isomorphism.
  \end{step}

  By \autoref{Prop_TameEstimate},
  it suffices to prove that $\SubellipticLaplacian_p$ is bijective.
  By \autoref{Prop_SubellipticEstimate+NeumannPoincareInequality}~\autoref{Prop_SubellipticEstimate+NeumannPoincareInequality_N},
  it is injective.
  For every $p \in \sU$ and $t \in (0,1]$,
  $\SubellipticLaplacian_{p,t} \co C_\perp^\infty(L) \to C_\perp^\infty(L)$ is elliptic.
  Since $\EnergyForm{p} \leq \EnergyForm{p,t}$,
  $\SubellipticLaplacian_{p,t}$ is injective.
  By the Riesz representation theorem and elliptic regularity,
  $\SubellipticLaplacian_{p,t}$ is surjective.
  By \autoref{Prop_TameEstimate},
  for every $k \in \N$,
  \begin{equation*}
    \norm{\SubellipticLaplacian_{p,t}^{-1} f}_{H^{k+\epsilon}} \lesssim_{k,p} \norm{f}_{H^k}.
  \end{equation*}
  For every $s,t \in (0,1]$,
  since
  \begin{equation*}
    \SubellipticLaplacian_{p,s}\paren{\SubellipticLaplacian_{p,t}^{-1} f - \SubellipticLaplacian_{p,s}^{-1}f}
    =
    (s-t)\Delta_g \SubellipticLaplacian_{p,t}^{-1} f,
  \end{equation*}
  by \autoref{Prop_TameEstimate},
  \begin{equation*}
    \norm{\SubellipticLaplacian_{p,t}^{-1} f - \SubellipticLaplacian_{p,s}^{-1}f}_{H^{k+\epsilon}}
    \lesssim_{k,p}
    \abs{s-t}
    \norm{f}_{H^{k+2}}.
  \end{equation*}
  Therefore,
  $\SubellipticLaplacian_{p,t}^{-1} f$ converges in $C_\perp^\infty(L)$ as $t \downarrow 0$
  and
  \begin{equation*}
    \SubellipticLaplacian_p \lim_{t \downarrow 0} \SubellipticLaplacian_{p,t}^{-1} f = f.
  \end{equation*}

  \begin{step}
    The inverses $(\SubellipticLaplacian_p^{-1})_{p \in \sU}$ form a smooth tame family.
  \end{step}

  By \cite[Part II, Theorem 3.1.1]{Hamilton1982:NashMoser},
  it suffices to prove that the map $\sU \times C_\perp^\infty(L) \to C_\perp^\infty(L), (p,f) \mapsto \SubellipticLaplacian_p^{-1}f$ is tame and continuous.
  \autoref{Prop_TameEstimate} shows that this map is tame,
  locally uniformly bounded in $p$, and
  continuous in $f$.
  This, together with the resolvent identity
  \begin{equation*}
    \SubellipticLaplacian_{p_0}^{-1} f_0 - \SubellipticLaplacian_p^{-1} f
    =
    \SubellipticLaplacian_p^{-1} \sqparen{f_0 - f + \paren{\SubellipticLaplacian_p - \SubellipticLaplacian_{p_0}}\SubellipticLaplacian_{p_0}^{-1}f_0},
  \end{equation*}
  shows that the map is jointly continuous.
\end{proof}


\section{Deformations of branched covers of special Lagrangians}
\label{Sec_DeformationsBranchedCoversOfSpecialLagrangians}

Throughout this section, assume the situation of \autoref{Thm_BranchedCoversOfSpecialLagrangians}.

Here is the computation announced in \autoref{Sec_Introduction}.
Set
\begin{equation*}
  \underline U \coloneq \set[\big]{ (t,\xi) \in \R \times T^*\Lambda : t\xi \in U }.
\end{equation*}

\begin{prop}
  \label{Prop_OmegaT}
  There is an $\underline \Omega \in \Omega^d(\underline U)$ such that:
  \begin{enumerate}
  \item
    \label{Prop_OmegaT_0}
    $\underline\Omega|_{\set{0} \times T^*\Lambda} = \Omega_g$.
  \item
    \label{Prop_OmegaT_T}
    For every $t \in \R^\times$,
    $\underline\Omega|_{\set{t} \times U_t} = \Omega_t$.
  \item
    \label{Prop_OmegaT_Exact+Primitive}
    For every $t \in \R$,
    \begin{equation*}
      \Omega_t ~\text{is exact} \qandq
      \Omega_t \wedge \omega = 0.
    \end{equation*}
  \end{enumerate}
\end{prop}

\begin{proof}
  The argument is based on expressing $\Omega_t$, defined in \autoref{Eq_OmegaT}, in a suitable form.
  Denote by $\xi\del_\xi \in \Vect(T^*\Lambda)$ the vector field generating the scaling action $\R^\times \circlearrowright T^*\Lambda$.
  Since $0^*\jmath^*\Im\theta = 0$,
  for every $t \in \R^\times$,
  \begin{align*}
    \Omega_t
    &=
      -t^{-1} \int_0^t \frac{\rd}{\rd s} s^*\jmath^*\Im\theta \,\rd s \\
    &=
      -\rd\sqparen*{t^{-1}\int_0^t s^{-1} s^*i_{\xi\del_\xi} \jmath^*\Im\theta \,\rd s}.
  \end{align*}
  Since
  \begin{equation*}
    \paren{s^{-1} s^*i_{\xi\del_\xi} \jmath^*\Im\theta}_\xi\paren{v_1,\ldots,v_{d-1}}
    =
    \paren{\jmath^*\Im\theta}_{s\xi}\paren{\xi, T_\xi s(v_1), \ldots, T_\xi s(v_{d-1})},
  \end{equation*}
  the integrand extends smoothly to $s = 0$.
  Moreover, at $s = 0$, it is given by
  \begin{equation*}
    \Im\theta_{\check \iota \circ \pi(\xi)}\paren{T_{0_{\pi(\xi)}}\jmath(\xi), T_\xi\paren{\check \iota \circ \pi}(v_1), \ldots, T_\xi\paren{\check \iota \circ \pi}(v_{d-1})}.
  \end{equation*}
  Since $[\check\iota]$ is special Lagrangian,
  only the normal component of $T_{0_{\pi(\xi)}}\jmath(\xi)$ contributes to this expression.
  The latter is given by $-JT_{\pi(\xi)}\check\iota\paren{\xi^\sharp}$.
  Since $\check\iota^*\Re\theta = \vol_g$ and $\theta$ is of type $(d,0)$,
  \begin{multline*}
    \Im\theta_{\check \iota \circ \pi(\xi)}\paren{T_{0_{\pi(\xi)}}\jmath(\xi), T_\xi\paren{\check \iota \circ \pi}(v_1), \ldots, T_\xi\paren{\check \iota \circ \pi}(v_{d-1})} \\
    =
    -\vol_g\paren{\xi^\sharp,T_\xi\pi(v_1), \ldots, T_\xi\pi(v_{d-1})}
    =
    -\Lambda_{g,\xi}\paren{v_1, \ldots,v_{d-1}}.
  \end{multline*}
  This proves \autoref{Prop_OmegaT_0}, \autoref{Prop_OmegaT_T}, and the exactness assertion in \autoref{Prop_OmegaT_Exact+Primitive}.
  Primitivity follows from $\Omega_g \wedge \omega = 0$ and $\jmath^*\theta \wedge \omega = 0$.
\end{proof}

Set $\sV
\coloneq
\set[\big]{
  (t,[\iota])
  \in \R \times \SpaceOfClosedNonDegenerateMultiValuedOneForms{\Lambda}
  :
  \im\iota \subseteq U_t
}$
and define $\DefiningSection \in \Gamma\paren{\sV, \pr_2^*\ExactFormBundle[d]}$ by
\begin{equation*}
  \DefiningSection(t,[\iota])
  \coloneq
  \paren{t,\sqparen{\iota;\iota^*\Omega_t}}.
\end{equation*}
By construction,
\begin{equation*}
  \DeformationNormalConeModuliSpaceOfSpecialLagrangians = \DefiningSection^{-1}(0).
\end{equation*}
The following analogue of \autoref{Prop_HodgeTypeOperator} is the final preparation for the proof of \autoref{Thm_BranchedCoversOfSpecialLagrangians}.

\begin{prop}
  \label{Prop_ThetaT}
  There is a $\underline \Theta \in \Omega^d(\underline U)$ with $\Theta_0=\Theta_g$ such that,
  for every
  $t \in \R$,
  $x \in U_t$,
  and $L \subseteq T_xU_t$ Lagrangian,
  the following hold:
  \begin{enumerate}
  \item
    \label{Prop_ThetaT_Calibrated=>Positive}
    If $\Omega_t|_L = 0$ and either $t\neq0$ or $\dim\paren{L\cap\ker T_x\pi}\leq2$,
    then $\Theta_t|_L \neq 0$ with $\Theta_t \coloneq \underline\Theta|_{\set{t} \times U_t}$.
  \item
    \label{Prop_ThetaT_Symmetric+NonNegativeDefinite}
    If $\Theta_t|_L \neq 0$,
    then $\HodgeTypeOperator{t}{L} \in \Hom\paren{L^*,\Wedge^{d-1} L^*}$ defined by
    \begin{equation*}
      \HodgeTypeOperator{t}{L}\alpha
      \coloneq
      \paren{i_v\Omega_t}|_L
      -
      \tfrac{\Omega_t|_L}{\Theta_t|_L}\paren{i_v\Theta_t}|_L
      \qwithq
      \paren{i_v\omega}|_L = \alpha
    \end{equation*}
    is symmetric and non-negative definite.
  \end{enumerate}
\end{prop}

\begin{proof}
  Decompose $TT^*\Lambda \iso \pi^*(T\Lambda\oplus T^*\Lambda)$ as in \autoref{Def_Theta} and,
  for every $t \in \R$,
  define
  \begin{equation*}
    \Theta_t
    \coloneq
    t^* \jmath^*\Re\theta
    -
    \frac12
    \left.\frac{\rd}{\rd s}\right|_{s=0}
    \Phi_s^*\Omega_t.
  \end{equation*}
  These assemble into $\underline \Theta \in \Omega^d(\underline U)$,
  by \autoref{Prop_OmegaT}.

  Since
  \begin{equation*}
    \Theta_0
    =
    \nu_{g,0}
    -
    \frac12
    \left.\frac{\rd}{\rd s}\right|_{s=0}\Phi_s^*\nu_{g,1}
    =
    \nu_{g,0}-\nu_{g,2}
    =
    \Theta_g,
  \end{equation*}
  by \autoref{Prop_Theta} and \autoref{Prop_PointwiseHodgeTypeOperator},
  the assertions hold for $t = 0$.
  It remains to prove the assertions for $t \in \R^\times$.

  \begin{proofsteps}
  \item[\autoref{Prop_ThetaT_Calibrated=>Positive}]
    Denote by $h$ the pullback of the Riemannian metric on $Z$ via $\jmath\circ t$.
    If $L$ is special Lagrangian, that is, $\Omega_t|_L = 0$,
    then
    $t^*\jmath^*\Re\theta|_L$ is a volume form of unit norm with respect to $h$.
    If $(e_1,\ldots,e_d)$ is a positive $h$--orthonormal basis of $L$,
    then,
    by direct computation with $(v_i,\xi_i) \coloneq e_i$,
    \begin{align*}
      \paren*{\left.\frac{\rd}{\rd s}\right|_{s=0} \Phi_s^*\Omega_t}(e_1,\ldots,e_d)
      &=
        \sum_{i=1}^d \Omega_t\paren[\big]{e_1,\ldots,e_{i-1},(\xi_i^\sharp,0),e_{i+1},\ldots,e_d} \\
      &=
        \sum_{i=1}^d \omega\paren{(\xi_i^\sharp,0),e_i}
        =
        -\sum_{i=1}^d \abs{\xi_i}_g^2.
    \end{align*}
    Therefore,
    $\Theta_t|_L$ is a positive multiple of $t^*\jmath^*\Re\theta|_L$.
  \item[\autoref{Prop_ThetaT_Symmetric+NonNegativeDefinite}]
    Since $\Theta_t \wedge \omega = 0$,
    $\HodgeTypeOperator{t}{L}$ is symmetric by the computation in the proof of \autoref{Prop_PointwiseHodgeTypeOperator}~\autoref{Prop_PointwiseHodgeTypeOperator_Symmetric}.

    Let $\alpha\in L^*$ with $\abs{\alpha}_h=1$.
    Choose a positive $h$--orthonormal basis $e_1,\ldots,e_d$ of $L$ with $e_1 = \alpha^\sharp$.
    Define $a+\iu b \in S^1 \subseteq \C$ by
    \begin{equation*}
      a+\iu b \coloneq \paren{t^*\jmath^*\theta}\paren{e_1,\ldots,e_d}
    \end{equation*}
    and set
    \begin{equation*}
      \tilde e_1\coloneq ae_1-bJ_te_1
      \qwithq
      J_t \coloneq t^*\jmath^*J.
    \end{equation*}

    The subspace $\tilde L = \Span{\tilde e_1,e_2,\ldots,e_d} \subseteq T_xU_t$ is special Lagrangian, because
    \begin{equation*}
      t^*\jmath^*\theta(\tilde e_1,e_2,\ldots,e_d)
      =(a-\iu b)(a+\iu b)=1.
    \end{equation*}
    Since $v=-tJ_te_1$ satisfies
    \begin{equation*}
      \paren{i_v\omega}|_L=\alpha
      \qandq
      \Omega_t(v,e_2,\ldots,e_d)=a,
    \end{equation*}
    by direct computation,
    \begin{align*}
      \frac{\alpha\wedge\HodgeTypeOperator{t}{L}\alpha}{\vol_h}
      &=
        a+\frac bt
        \frac{\Theta_t(v,e_2,\ldots,e_d)}
        {\Theta_t(e_1,\ldots,e_d)}
      \\
      &=
        \frac{\Theta_t(\tilde e_1,e_2,\ldots,e_d)}
        {\Theta_t(e_1,\ldots,e_d)}.
    \end{align*}
    The latter is positive by \autoref{Prop_ThetaT_Calibrated=>Positive}.
    \qedhere
  \end{proofsteps}
\end{proof}

\begin{proof}[Proof of \autoref{Thm_BranchedCoversOfSpecialLagrangians}]
  Over $\R^\times$,
  the assertion follows from \cite[Theorem 3.6]{McLean1998} and \cite[§4]{Hitchin1997:ModuliSpaceOfSpecialLagrangianSubmanifolds}.

  Set
  \begin{equation*}
    \PositiveLocus
    \coloneq
    \set[\big]{
      (t,[\iota]) \in \sV
      :
      \iota^*\Theta_t ~\text{is a volume form}
    }
  \end{equation*}
  and define the covariant derivative $\nabla$ on $\pr_2^*\ExactFormBundle|_{\PositiveLocus}$ as in \autoref{Def_PositiveLocus} with $\Theta_t$ instead of $\Theta_g$.

  \autoref{Prop_CovariantDerivativeOfDefiningSection} carries over to the present situation.
  By \autoref{Prop_CoefficientModuleSatisfiesHormandersCondition},
  \autoref{Prop_ThetaT},
  \autoref{Thm_TameFamilyOfInversesForSubellipticLaplacians}, and
  the argument in the proof of \autoref{Prop_ExtendedCovariantDerivativeIsIsomorphism},
  there is an open neighbourhood $\set{0} \times \SpaceOfHarmonicNonDegenerateMultiValuedOneForms{\Lambda,g} \subseteq \sU \subseteq \PositiveLocus$ such that
  \begin{equation*}
    \paren{\nabla\DefiningSection ~ \GaussManinConnection\PeriodSection}
    \co \pr_2^*T\SpaceOfClosedNonDegenerateMultiValuedOneForms{\Lambda}|_\sU
    \to
    \pr_2^*\paren{\ExactFormBundle[d] \oplus \DeRhamCohomologyBundle[1]}|_\sU
  \end{equation*}
  is an isomorphism of tame Fréchet space orbibundles.
  This, together with the Nash--Moser submersion theorem \cite[Theorem 4.7]{HeParkerWalpuski2026}, implies the (remaining) assertion.
\end{proof}


\appendix
\section{The relation to harmonic \texorpdfstring{$W_\loc^{1,2}$ multi-valued $1$}{W1,2loc multi-valued 1}--forms}
\label{Sec_RelationHarmonicW12MultiValued1Forms}

This appendix discusses the relation between harmonic non-degenerate multi-valued $1$--forms and harmonic $W_\loc^{1,2}$ multi-valued $1$--forms in the sense of \citet{Almgren2000:BigRegularityPaper}.
Throughout, assume that $X$ is connected.

\begin{prop}
  \label{Prop_HarmonicNonDegenerateMultiValuedOneForm=>Almgren+ZModTwoZ}
  For every $(g,[\iota \co L \imm T^*X]) \in \SpaceOfHarmonicNonDegenerateMultiValuedOneForms{X}$ with multiplicity $k = k[\iota] \in \N$,
  there is a unique
  \begin{equation*}
    \alpha = \alpha([\iota]) \in W_\loc^{1,2}\Gamma\paren{X,\Sym^k T^*X}
    \qwithq
    \rd\alpha = 0
    \qandq
    \rd^*\alpha = 0
  \end{equation*}
  such that, for almost every $x \in X$,
  \begin{equation}
    \label{Eq_IotaToAlpha}
    \alpha(x) = \sum_{y \in \paren{\pi \circ \iota}^{-1}(x)} \Current{\iota(y)}.
  \end{equation}
  Moreover, the barycentre $\overline\alpha$ is smooth and harmonic, that is, $\overline\alpha \in \sH^1(X,g)$.
\end{prop}

\begin{proof}
  By \autoref{Prop_StructureOfProjectionMap}~\autoref{Prop_StructureOfProjectionMap_RamificationIndex},
  \autoref{Eq_IotaToAlpha} defines a $k$--valued $1$--form $\alpha$ on $X \setminus \BranchingLocus{\iota}$
  with harmonic local branches $\alpha_1,\ldots,\alpha_k$.
  
  Consider an open subset $U \subseteq X$ whose closure is compact and contained within a coordinate neighbourhood.
  Since $\pi \circ \iota$ is proper,
  $\alpha$ is bounded on $U \setminus \BranchingLocus{\iota}$.
  Moreover, by \autoref{Eq_OmegaThetaViaS},
  \begin{equation*}
    \int_{U \setminus \BranchingLocus{\iota}} \sum_{i=1}^k \abs{\nabla\alpha_i}^2 \, \vol_g
    \leq
    2 \int_{\paren{\pi \circ \iota}^{-1}(U)} \iota^*\Theta_g
    < \infty.
  \end{equation*}
  
  By \autoref{Prop_StructureOfProjectionMap}~\autoref{Prop_StructureOfProjectionMap_NoFolds} and Sard's theorem applied to the compositions of $\pi \circ \iota$ with the coordinate projections $p_j \co U \to \bR^{d-1}$,
  for almost every $y \in \R^{d-1}$,
  the fibre $p_j^{-1}(y)$ is disjoint from $\BranchingLocus{\iota}$.
  In the coordinate trivialisation,
  for every $T \in \Sym^k \R^d$,
  the distance $d(\alpha,T)$ 
  is locally Lipschitz along such a fibre and satisfies
   \begin{equation*}
     \abs{\del_j d(\alpha,T)}^2 
     \leq \sum_{i=1}^k\abs{\del_j\alpha_i}^2
     \lesssim_U 1 + \sum_{i=1}^k\abs{\nabla\alpha_i}^2.
   \end{equation*}
   Integration by parts along the fibres and Fubini’s theorem identify these derivatives with the weak derivatives on $U$.
   Therefore,
   $\alpha \in W_\loc^{1,2}\Gamma\paren{X,\Sym^k T^*X}$,
   by \cite[Definition 0.5]{DeLellisSpadaro2011:QValuedRevisited}.
   Uniqueness and harmonicity follow because $\BranchingLocus{\iota}$ has measure zero.
   
   Since the barycentre map is Lipschitz,
   $\overline\alpha \in W_\loc^{1,2}\Omega^1(X)$ and is weakly harmonic because the local branches are.
   By elliptic regularity,
   $\overline\alpha$ is smooth.
\end{proof}


\section{The relation to harmonic non-degenerate \texorpdfstring{$\ZModTwoZ$ $1$}{Z/2Z 1}--forms}
\label{Sec_RelationHarmonicZModTwoZ1Forms}

This appendix discusses the relation between harmonic non-degenerate $2$--valued $1$--forms and harmonic non-degenerate $\Z/2\Z$ $1$--forms.
Throughout, let $g \in \SpaceOfMetrics[X]$.

\begin{definition}[cf.~{\cites[§1]{Donaldson2021}[Definitions~7 and~8]{He2022:DeformBranchedSLags}}]
  \label{Def_HarmonicZModTwoZOneForm}
  A \defined{harmonic non-degenerate $\ZModTwoZ$ $1$--form} on $(X,g)$ is a pair $(\fl,\alpha)$ consisting of:
  \begin{enumerate}
  \item
    a Euclidean line bundle $\fl$ over $X \setminus \BranchingLocus{\fl}$ where $\BranchingLocus{\fl} \subseteq X$ is a closed submanifold of codimension two,
    and
  \item
    a twisted harmonic $1$--form $\alpha \in \Omega^1\paren{X \setminus \BranchingLocus{\fl}; \fl}$ satisfying
    \begin{equation*}
      \abs{\alpha} \asymp \dist[g]{\BranchingLocus{\fl}}{\cdot}^{1/2}
      \quad
      \text{locally near}
      \quad
      \BranchingLocus{\fl}.
      \qedhere
    \end{equation*}
  \end{enumerate}
\end{definition}

Using \autoref{Prop_StructureOfProjectionMap},
it is relatively straightforward to see that every harmonic non-degenerate $2$--valued $1$--form with zero barycentre induces a harmonic non-degenerate $\ZModTwoZ$ $1$--form.
The converse also holds,
as is explained in the following.

\begin{prop}
  \label{Prop_RelationHarmonicZModTwoZ1Forms}
  For every harmonic non-degenerate $\ZModTwoZ$ $1$--form $(\fl,\alpha)$ on $(X,g)$,
  the double cover $\varpi^\circ \co \tilde X^\circ \coloneq \set{ v \in \fl : \abs{v} = 1} \to X \setminus \BranchingLocus{\fl}$ and
  the map $\iota^\circ \co \tilde X^\circ \to T^*X$ induced by $\alpha$
  extend to
  a smooth branched double cover $\varpi \co \tilde X \to X$ and
  an immersion $\iota \co \tilde X \imm T^*X$ such that $\pi \circ \iota = \varpi$ and $[\iota]$ is a harmonic non-degenerate $2$--valued $1$--form with zero barycentre.
\end{prop}

The proof relies on the following regularity result.

\begin{lemma}
  \label{Lem_SmoothAfterPullback}
  Let $\BranchingLocus{\fl} \subseteq X$ be a closed submanifold of codimension two, and $\fl$ a Euclidean line bundle over $X \setminus \BranchingLocus{\fl}$.
  The double cover $\varpi^\circ \co \tilde X^\circ \to X \setminus \BranchingLocus{\fl}$ extends to a smooth branched double cover $\varpi \co \tilde X \to X$ such that the following holds:
  if $E$ is a Euclidean vector bundle over $X$ equipped with an orthogonal covariant derivative $\nabla$ and $\sR \in \Gamma\paren{X,\End(E)}$, and
  $s \in \Gamma\paren{X \setminus \BranchingLocus{\fl},E \otimes \fl}$ satisfies
  \begin{equation*}
    \nabla^*\nabla s + \sR s = 0
    \qandq
    \abs{s} \lesssim \dist[g]{\BranchingLocus{\fl}}{\cdot}^{1/2},
  \end{equation*}
  then $\paren{\varpi^\circ}^*s \in \Gamma\paren{\tilde X^\circ,\varpi^*E}$ extends smoothly to $\tilde X$.
\end{lemma}

\begin{proof}[Proof of \autoref{Prop_RelationHarmonicZModTwoZ1Forms} assuming \autoref{Lem_SmoothAfterPullback}]
  By Bochner's formula,
  \begin{equation*}
    \nabla^*\nabla\alpha + \Ric(\alpha) = 0.
  \end{equation*}
  Therefore, \autoref{Lem_SmoothAfterPullback}
  equips $\tilde X$ with a smooth structure for which $\varpi$ is a smooth branched cover and $\iota^\circ$ extends to a smooth map $\iota \co \tilde X \to T^*X$ satisfying $\pi \circ \iota = \varpi$.

  By the growth condition on $\alpha$,
  $\iota$ is an immersion.
  Since $\alpha$ is harmonic,
  \begin{equation*}
    \iota^*\omega = 0 \qandq \iota^*\Omega_g = 0.
  \end{equation*}
  Since $\varpi$ is a smooth branched double cover,
  $[\iota]$ is a harmonic non-degenerate $2$--valued $1$--form.
\end{proof}

\begin{proof}[Proof of \autoref{Lem_SmoothAfterPullback}]
  By elliptic regularity,
  $\tilde s^\circ \coloneq \paren{\varpi^\circ}^*s$ extends across every component of $\BranchingLocus{\fl}$ around which $\fl$ has trivial monodromy.
  Therefore,
  it suffices to assume that $\fl$ has monodromy $-1$ around every component of $\BranchingLocus{\fl}$.
  The proof proceeds in three steps.

  \begin{step}
    \label{Step_Lem_SmoothAfterPullback_SmoothStructure}
    Construction of the smooth structure.
  \end{step}

  \cite[Proof of Lemma 7.18]{HambletonHausmann2011} constructs a smooth structure on $\tilde X$ with respect to which $\varpi^\circ$ extends to $\varpi$.
  Furthermore,
  if $(t,\zeta)$ with $t = (t_1,\ldots,t_{d-2})$ and $\zeta = \zeta_1 + \iu \zeta_2$ are Fermi coordinates with respect to $g$ and $\BranchingLocus{\fl}$ on an open subset $U \subseteq X$,
  then there are corresponding coordinates $(t,z)$ with $z = z_1 + \iu z_2$ on $\tilde U \coloneq \varpi^{-1}(U)$ with respect to which
  \begin{equation*}
    \varpi(t,z) = (t,z^2).
  \end{equation*}

  \begin{step}
    \label{Step_Lem_SmoothAfterPullback_Hypoelliptic}
    A hypoelliptic operator.
  \end{step}

  By hypothesis,
  \begin{equation*}
    P^\circ \tilde s^\circ = 0
    \qwithq
    P^\circ
    \coloneq
    \abs{z}^2\paren{\varpi^\circ}^*\paren{\nabla^*\nabla + \sR}.
  \end{equation*}
  Evidently,
  $\abs{z}^2\paren{\varpi^\circ}^*\sR$ extends to $\tilde X$.
  The following discussion shows that the same is true for $P^\circ$.

  The rough Laplacian is $\nabla^*\nabla = -\tr_g\nabla^{(2)} = -\inner{g^{-1}}{\nabla^{(2)}}$.
  By the generalised Gauss lemma \cites[Lemma 2.11]{Gray2004}[Chapter X, §IV, Theorem II]{Cartan1983},
  the inverse metric has the form~
  \begin{align*}
    g^{-1}
    &=
      \del_{\zeta_1}^{\otimes 2} + \del_{\zeta_2}^{\otimes 2}
      + a \del_\alpha^{\otimes 2}
      + \sum_{i=1}^{d-2}\sum_{j=1}^2 b_{ij}\del_{t_i} \odot \del_{\zeta_j}
      + \sum_{i,j=1}^{d-2} c_{ij}\del_{t_i} \otimes \del_{t_j}
  \end{align*}
  with $\del_\alpha \coloneq \zeta_1 \del_{\zeta_2} - \zeta_2 \del_{\zeta_1}$ and smooth coefficients $a$, $(b_{ij})$, and $(c_{ij})$.
  Here $u \odot v \coloneq \tfrac12(u \otimes v + v \otimes u)$.
  Therefore,
  $\abs{z}^2\varpi^* g^{-1}$ admits a representation as
  \begin{equation*}
    \abs{z}^2\varpi^* g^{-1} = \sum_{i=1}^N v_i \otimes v_i
  \end{equation*}
  with $v_1,\ldots,v_N \in \Vect(\tilde U)$ generating the same $C^\infty$--submodule as $z_i\del_{t_j}$, $\del_{z_i}$ ($i \in \set{1,2}$, $j \in \set{1,\ldots,d-2}$).
  This implies that $P^\circ$ extends to a differential operator $P$ on $\tilde X$.
  By Hörmander's theorem \cites[Theorem 1.1]{Hormander1967:Hypoelliptic}[§2.2.1]{Ponge2008:HeisenbergCalculus},
  $P$ is hypoelliptic.

  \begin{step}
    \label{Step_Lem_SmoothAfterPullback_Distributional}
    $\tilde s^\circ$ extends smoothly to $\tilde X$.
  \end{step}

  By assumption, $\tilde s^\circ$ extends continuously to $\tilde X$, and
  its extension $\tilde s$ satisfies $\abs{\tilde s} \lesssim_{\tilde U} \abs{z}$.
  For every $\epsilon \in (0,1]$,
  choose a cutoff function $\chi_\epsilon$ which vanishes for $\abs{z} \leq \epsilon$, equals one for $\abs{z} \geq 2\epsilon$, and satisfies $\abs{\nabla^{(j)}\chi_\epsilon} \lesssim \epsilon^{-j}$ for $j \in \set{1,2}$.
  Denote by $P^*$ the formal adjoint with respect to $\rd t\,\rd z$.
  For every $\varphi \in \Gamma_c\paren{\tilde U,\varpi^*E}$,
  \begin{equation*}
    \abs*{\int_{\tilde U}\chi_\epsilon\inner{\tilde s}{P^*\varphi}\,\rd t\,\rd z}
    =
    \abs*{\int_{\tilde U}\inner{\tilde s}{\LieBracket{P^*}{\chi_\epsilon}\varphi}\,\rd t\,\rd z}
    \lesssim
    \epsilon.
  \end{equation*}
  Therefore,
  $P\tilde s = 0$ in the sense of distributions.
  Since $P$ is hypoelliptic,
  $\tilde s$ is smooth.
\end{proof}


\printbibliography

\end{document}
